\documentclass[final]{amsart}

\usepackage{amssymb}
\usepackage[foot]{amsaddr}
\usepackage[utf8]{inputenc}
\usepackage[T1]{fontenc}
\usepackage[english, ]{babel}
\usepackage{csquotes} \usepackage{xargs}
\usepackage{graphicx, color}
\usepackage[figurewithin=none]{caption,subcaption} \graphicspath{{./}} \usepackage{enumitem}
\setlist{itemsep=3pt,topsep=3pt}

\usepackage{array, multirow}

\newcommandx{\mr}[3][2=*]{\multirow{#1}{#2}{#3}}

\usepackage{tabularx}
\newcolumntype{C}{>{\centering\arraybackslash}X}\newcolumntype{R}{>{\raggedleft\arraybackslash}X}\newcolumntype{L}{>{\raggedright\arraybackslash}X}\newcolumntype{M}{>{$}C<{$}}\newcolumntype{\ll}[1]{>{\hsize=#1\hsize\raggedright\arraybackslash}X}\newcolumntype{\rr}[1]{>{\hsize=#1\hsize\raggedleft\arraybackslash}X}\newcolumntype{\cc}[1]{>{\hsize=#1\hsize\centering\arraybackslash}X}\newcolumntype{\mm}[1]{>{$}\cc{#1}<{$}}

\usepackage[notcite,notref]{showkeys} \usepackage{hyperref}
\hypersetup{
  hidelinks,
  colorlinks  = true,    urlcolor    = blue,    linkcolor   = blue,    citecolor   = blue,      pdfauthor   = {Antonio Montero and Primoz Potocnik},pdfsubject  = {Fixity, permutation groups, vertex-transitive graphs},pdftitle    = {Small motion},  }
  
\usepackage[capitalise,noabbrev, nameinlink, poorman]{cleveref}
\usepackage{tikz-cd}
  \tikzcdset{column sep=tiny,
  row sep = tiny,
  }
  
\usepackage[loadshadowlibrary,shadow, colorinlistoftodos,textsize=tiny]{todonotes}
  \setuptodonotes{fancyline, backgroundcolor=gray!10,bordercolor=blue}

\newcommandx{\inline}[2][1=]{\todo[inline, #1]{#2}}

  \makeatletter
  \providecommand\@dotsep{5}
  \makeatother

\theoremstyle{plain}
\newtheorem{thm}{Theorem}[section]
\newtheorem*{thm*}{Theorem}

\newtheorem{lem}[thm]{Lemma}
\newtheorem{lemma}[thm]{Lemma}
\newtheorem{prop}[thm]{Proposition}

\newtheorem{coro}[thm]{Corollary}

\theoremstyle{definition}

\newtheorem{definition}[thm]{Definition}

\newtheorem{example}[thm]{Example}
\newtheorem{problem}{Problem}
\newtheorem{prob}[problem]{Problem}

\theoremstyle{remark}

\newtheorem{rem}[thm]{Remark}
\newtheorem{remark}[thm]{Remark}
\newtheorem{claim}{Claim}

\numberwithin{equation}{section}

\Crefname{thm}{Theorem}{Theorems}
\Crefname{theorem}{Theorem}{Theorems}
\Crefname{lem}{Lemma}{Lemmas}
\Crefname{lemma}{Lemma}{Lemmas}
\Crefname{prop}{Proposition}{Propositions}
\Crefname{proposition}{Proposition}{Propositions}
\Crefname{coro}{Corollary}{Corollaries}
\Crefname{corollary}{Corollary}{Corollaries}
\Crefname{defn}{Definition}{Definitions}
\Crefname{definition}{Definition}{Definitions}
\Crefname{exam}{Example}{Examples}
\Crefname{example}{Example}{Examples}
\Crefname{note}{Note}{Notes}
\Crefname{rem}{Remark}{Remarks}
\Crefname{remark}{Remark}{Remarks}
\Crefname{prob}{Problem}{Problems}
\Crefname{problem}{Problem}{Problems}
\Crefname{claim}{Claim}{Claims}

\usepackage{xargs}

\renewcommand{\leq}{\leqslant} \renewcommand{\geq}{\geqslant}
\renewcommand{\epsilon}{\varepsilon} \renewcommand{\theta}{\vartheta} \renewcommand{\subset}{\subseteq}  
\renewcommand{\{}{\lbrace}
\renewcommand{\}}{\rbrace}
\newcommand{\sm}{\setminus}
\renewcommand{\bar}{\overline}

\let\Wr\wr
\renewcommandx{\wr}{\mathbin{\mathrm{wr}}}

\newcommand{\bF}{\mathbb{F}}
\newcommand{\bN}{\mathbb{N}}

\newcommand{\bZ}{\mathbb{Z}}

\newcommand{\cB}{\mathcal{B}}

\newcommand{\cD}{\mathcal{D}}

\newcommand{\cM}{\mathcal{M}}

\newcommand{\cP}{\mathcal{P}}

\newcommand{\cyvec}[1]{\bar{\mathrm{#1}}}

\newcommand{\vid}{\cyvec{1}}

\newcommand{\rr}{\varrho}

\newcommand{\bfg}{\mathbf{g}}

\DeclareMathOperator{\supp}{supp}
\DeclareMathOperator{\sgn}{sgn}

\DeclareMathOperator{\fun}{Fun}
\DeclareMathOperator{\aut}{Aut}

\DeclareMathOperator{\mindeg}{\mu}

\DeclareMathOperator{\iinf}{Inf}

\DeclareMathOperator{\sym}{Sym}
\DeclareMathOperator{\alt}{Alt}
\DeclareMathOperator{\dih}{Dih}

\DeclareMathOperator{\agl}{AGL}
\DeclareMathOperator{\aGl}{A\Gamma L}
\DeclareMathOperator{\pgl}{PGL}
\DeclareMathOperator{\psl}{PSL}
\DeclareMathOperator{\pGl}{P\Gamma L}
\DeclareMathOperator{\gl}{GL}

\newcommandx{\dtwoSM}[2][1=\cM, 2=s] {\hat{2}#2^{#1 - 1}}

\newcommandx{\cube}[1][1=m, usedefault]{C_2^{#1} \rtimes \sym(#1)}
\newcommandx{\evencube}[1][1=m, usedefault]{\left(C_{2}^{#1} \rtimes \sym(#1) \right)^{+}}
\newcommandx{\cubeWr}[1][1=m, usedefault]{C_{2} \wr \sym(#1)}
\newcommandx{\evencubeWr}[1][1=m, usedefault]{\left(C_{2} \wr \sym(#1)\right)^{+}}
\newcommand{\ab}[1][m]{{C_{2}^{#1}}}
\newcommand{\evenAb}[1][m]{{(C_{2}^{#1})^{+}}}

\newcommandx{\impWr}[3][1=G, 2=, 3=G^{\cB}, usedefault]{#1_{#2}^{#2} \wr {#3}}

\newcommand{\lex}[2]{#1 \Wr #2}

\newcommandx{\Inf}[5][1=\lambda, 2=\kappa, 3=\Sigma, 4=\cP,, 5=m, usedefault] {\iinf^{#1}_{#2}(#3,#4,#5)}

\newcommand{\pleq}{\leq}

\newcommand{\sflip}{\tau}

\usepackage[sorting=nyt, backend=bibtex, style=numeric,defernumbers=true,url=false,eprint=false,giveninits=true,
maxbibnames=9,isbn=false,sortcites = true,citestyle=numeric-comp,
]{biblatex}
\begin{document}

\title[small motion and small minimal degree]{Vertex-transitive graphs with small motion and transitive permutation groups with small minimal degree}

\author{Antonio Montero $^{a}$}
\address{$^a$ Faculty of Mathematics and Physics, University of Ljubljana, SI-1000 Ljubljana, Slovenia}

\author{Primož Potočnik $^{a,b}$}
\address{$^b$ Institute of Mathematics, Physics and Mechanics, Jadranska 19, SI-1000
Ljubljana, Slovenia}

\address{Research supported in part by the Slovenian Research and Inovation Agency (ARIS) program no. P1-0294 and project no. N1-0216 }
\email{antonio.montero@fmf.uni-lj.si, primoz.potocnik@fmf.uni-lj.si}


\begin{abstract}
  The motion of a graph is the minimum number of vertices that are moved by a non-trivial automorphism. Equivalently, it can be defined as the minimal degree of its automorphism group (as a permutation group on the vertices). In this paper we develop some results on permutation groups (primitive and imprimitive) with small minimal degree. As a consequence of such results we classify vertex-transitive graphs whose motion is $4$ or a prime number.
\end{abstract}

\maketitle

\section{Introduction and main results} \label{sec:intro}
The \emph{minimal degree} $\mindeg(G)$ (sometimes called the \emph{motion}) of a permutation group  $G\leq \sym(\Omega)$ is defined to be the smallest number of points moved by any
non-identity element of $G$, that is
\[\mindeg(G) = \min\left\{ |\supp(x)| : x \in G \sm \left\{ 1 \right\}  \right\} \]
where $\supp(x) = \left\{ \omega \in \Omega : \omega^x \neq \omega \right\} $ is the \emph{support} of $x$.
Recently, we have witnessed a growing interest in the minimal degree of permutation groups that appear as automorphism groups of graphs; see \cite{Babai_2014_AutomorphismGroupsStrongly,Babai_2015_AutomorphismGroupsStrongly,PotocnikSpiga_2021_NumberFixedPoints,LehnerPotocnikSpiga_2021_FixityArcTransitive,BarbieriGrazianSpiga_2023_NumberFixedEdges,BarbieriPotocnik_2024_VertexPrimitiveDigraphs}.
Following the terminology of \cite{RussellSundaram_1998_NoteAsymptoticsComputational,ConderTucker_2011_MotionDistinguishingNumber}, we define the \emph{motion} $\mu(\Gamma)$ of a graph $\Gamma$ as the minimal degree of its automorphism group, that is \[\mu(\Gamma) = \mindeg(\aut(\Gamma)),\]
where $\aut(\Gamma)$ is viewed as a permutation group on the vertex-set of $\Gamma$.
The main goal of this paper is to investigate vertex-transitive graphs of small motion; see \cref{thm:graphs}.
If $\Gamma_{1}$ and $\Gamma_{2}$ are graphs, we let $\lex{\Gamma_{2}}{\Gamma_{1}}$ (also often denoted $\Gamma_{1}[\Gamma_{2}]$) and $\Gamma_{1} \square \Gamma_{2}$ denote the lexicographic and the cartesian product of $\Gamma_{1}$ and $\Gamma_{2}$, respectively (see \cite{ImrichKlavzar_2000_ProductGraphsStructure}).

\begin{thm}\label{thm:graphs}
  Let $\Gamma$ be a vertex-transitive graph on $n$ vertices. 
Then:
  \begin{enumerate}
    \item \label{item:graphs_mindeg2_} 
    If $\aut(\Gamma)$ contains a $2$-cycle (or equivalently if $\mu(\Gamma) = 2$) then
$\Gamma \cong \lex{K_{m}}{\Theta}$ or $\Gamma \cong \lex{(mK_{1})}{\Theta}$ with
    $\Theta$ a vertex-transitive graph on $k$ vertices and $m \geq 2$ such that $n =mk$.  
    Conversely, for every vertex-transitive graph $\Theta$ and $m \geq 2$ the graphs $\lex{K_{m}}{\Theta}$ and $\lex{(mK_{1})}{\Theta}$ have motion $2$.
    \item \label{item:graphs_mindegp_} If $\aut(\Gamma)$ contains a $p$-cycle for some prime number $p \geq 3$, then either:
    \begin{enumerate}
      \item $\mu(\Gamma) = 2$ and thus $\Gamma$ is one of the graphs described in (\ref{item:graphs_mindeg2}); or
      \item $\Gamma$ is isomorphic to $\lex{\Sigma_{p}}{\Theta}$ with $\Sigma_{p}$ a circulant graph with $p$ vertices and $\Theta$ is a vertex-transitive graph; in this case $\mu(\Gamma) = 2$ if $\Sigma_{p}$ is isomorphic to $ K_{p} $ or $pK_{1}$, and $\mu(\Gamma) = p-1$ otherwise. 
    \end{enumerate}
\item \label{item:graphs_mindeg4_} If $\mu(\Gamma) = 4$ then one of the following holds:
      \begin{enumerate}
        \item 
        \label{item:graphs_mindeg4_prim_}
        $\Gamma \cong \lex{C_{5}}{\Theta}$ with $\Theta$ a vertex-transitive graph; or \item \label{item:graphs_mindeg4_cross1_}
         $\Gamma \cong \lex{(K_{m} \square K_{2})}{\Theta}$ or $\Gamma \cong \lex{\bar{K_{m} \square K_{2}}}{\Theta}$ with $m\geq 3  $ and $\Theta$ a vertex-transitive graph; or
\item \label{item:graphs_mindeg4_cross2_} $\Gamma \cong  \Inf$, defined in \cref{def:inf} below, for some 
        $m \geq 2$, $\lambda,\kappa \in \bZ_{2}$, a graph $\Sigma$, and a partition $\cP$ of $V\Sigma$ with blocks of size $2$ satisfying some natural conditions precised in \cref{lem:inf}.
        \end{enumerate}
      Conversely, for every vertex-transitive graph $\Theta$, every vertex-transitive graph $\Sigma$ admitting a $\aut(\Sigma)$-invariant partition $\cP$ of its vertices with blocks of size $2$ satisfying the conditions in \cref{lem:inf}, and  $m \geq 2$, $\lambda, \kappa \in \bZ_{2}$, the graphs $\lex{C_5}{\Theta}$, $\lex{(K_{m} \square K_2)}{\Theta}$, $\lex{\bar{K_{m} \square K_{2}}}{\Theta}$ and $\Inf$  have motion $4$.
  \end{enumerate}
\end{thm}

\begin{coro}
  The motion of a vertex-transitive graph is never an odd prime.
\end{coro}

\begin{definition}
  \label{def:inf}
  Let $\Sigma$ be a graph with vertex set $\Omega$ and let $\cP$ be a partition of $\Omega$ into sets of size $2$.
  Let $m \geq 2$ be an integer and let $\lambda,\kappa \in \bZ_{2} $.
  The graph $\Inf$ is the graph whose vertex set is $\Omega \times \left\{ 1, \dots, m \right\} $ and the edges are given by
  \[
    (\alpha, i) \sim (\beta, j) \iff
    \begin{cases}
      \left\{ \alpha, \beta \right\} \not\in \cP
      \text{ and }
      \alpha \sim \beta \text{ in } \Sigma
      \\
      \left\{ \alpha,\beta \right\} \in \cP,
      \lambda = 1,
      \text{ and }
      i=j, \\
      \left\{ \alpha,\beta \right\} \in \cP,
      \lambda = 0,
      \text{ and }
      i\neq j,\\
      \alpha = \beta \text{ and } \kappa = 1.
    \end{cases}
    \]
  \end{definition} 

\begin{figure}
  \begin{subfigure}[b]{0.45\textwidth}
    \centering
    \includegraphics[width=\textwidth]{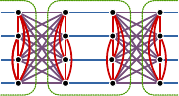}
    \caption{$\lambda=1, \kappa=1$}
  \end{subfigure}
  \hfill
  \begin{subfigure}[b]{0.45\textwidth}
    \centering
    \includegraphics[width=\textwidth]{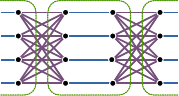}
    \caption{$\lambda=1, \kappa=0$}
\end{subfigure}

  \begin{subfigure}[b]{0.45\textwidth}
    \centering
    \includegraphics[width=\textwidth]{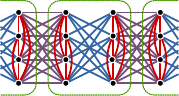}
    \caption{$\lambda=0, \kappa=1$}
\end{subfigure}
  \hfill
  \begin{subfigure}[b]{0.45\textwidth}
    \centering
    \includegraphics[width=\textwidth]{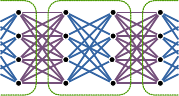}
    \caption{$\lambda=0, \kappa=0$}
\end{subfigure}
  \caption{Local view of the graph $\Inf$ with $m=4$.}
  \label{fig:inf}
\end{figure}

 \begin{remark}
  Informally speaking, to obtain $\Inf$ every vertex of $\Sigma$ is turned into an independent set of $m$ vertices (if $\kappa =0$) or into a complete graph $K_{m}$ (if $\kappa = 1$). This set of $m$ vertices of $\Inf$ resulting from a single vertex $v$ of $\Sigma$ will informally be called a {\em fibre above} $v$.
  The edges in $\Sigma$ that do not define an element of $\cP$ are turned into complete bipartite graphs, just as in the lexicographic product.
  The pairs $\{u,v\}\in \cP$ induce a matching between the fibres above $u$ and $v$ (if $\lambda=1$) or complete bipartite graphs without a matching (if $\lambda= 0$), regardless of whether the vertices $u$ and $v$ are adjacent in $\Sigma$; see \cref{fig:inf}.
Note that if $\Sigma'$ is obtained from $\Sigma$ by removing an edge forming a pair in $\cP$, then $\Inf \cong \Inf[][][\Sigma']$. Hence we could assume that $\alpha \sim \beta$ whenever $\left\{ \alpha, \beta \right\} \in \cP$. 
  Observe also that
  \[\bar{\Inf} =
  \Inf[\lambda+1][\kappa+1][\bar{\Sigma}].
  \]
\end{remark} 

\begin{figure}
  \centering
  \begin{subfigure}{\textwidth}
    \centering
    \begingroup \makeatletter \providecommand\color[2][]{\errmessage{(Inkscape) Color is used for the text in Inkscape, but the package 'color.sty' is not loaded}\renewcommand\color[2][]{}}\providecommand\transparent[1]{\errmessage{(Inkscape) Transparency is used (non-zero) for the text in Inkscape, but the package 'transparent.sty' is not loaded}\renewcommand\transparent[1]{}}\providecommand\rotatebox[2]{#2}\newcommand*\fsize{\dimexpr\f@size pt\relax}\newcommand*\lineheight[1]{\fontsize{\fsize}{#1\fsize}\selectfont}\ifx\svgwidth\undefined \setlength{\unitlength}{174.33070866bp}\ifx\svgscale\undefined \relax \else \setlength{\unitlength}{\unitlength * \real{\svgscale}}\fi \else \setlength{\unitlength}{\svgwidth}\fi \global\let\svgwidth\undefined \global\let\svgscale\undefined \makeatother \begin{picture}(1,0.36105592)\lineheight{1}\setlength\tabcolsep{0pt}\put(0,0){\includegraphics[width=\unitlength,page=1]{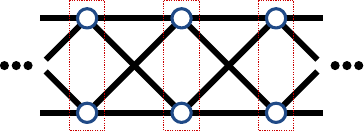}}\end{picture}\endgroup      \caption{$\mathrm{PX}(r,1) \cong \lex{(2 K_{1})}{C_{r}}$}
  \end{subfigure}
  
  \begin{subfigure}{\textwidth}
    \centering
    \begingroup \makeatletter \providecommand\color[2][]{\errmessage{(Inkscape) Color is used for the text in Inkscape, but the package 'color.sty' is not loaded}\renewcommand\color[2][]{}}\providecommand\transparent[1]{\errmessage{(Inkscape) Transparency is used (non-zero) for the text in Inkscape, but the package 'transparent.sty' is not loaded}\renewcommand\transparent[1]{}}\providecommand\rotatebox[2]{#2}\newcommand*\fsize{\dimexpr\f@size pt\relax}\newcommand*\lineheight[1]{\fontsize{\fsize}{#1\fsize}\selectfont}\ifx\svgwidth\undefined \setlength{\unitlength}{219.68503937bp}\ifx\svgscale\undefined \relax \else \setlength{\unitlength}{\unitlength * \real{\svgscale}}\fi \else \setlength{\unitlength}{\svgwidth}\fi \global\let\svgwidth\undefined \global\let\svgscale\undefined \makeatother \begin{picture}(1,0.31232182)\lineheight{1}\setlength\tabcolsep{0pt}\put(0,0){\includegraphics[width=\unitlength,page=1]{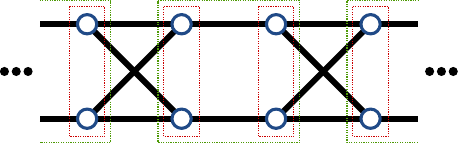}}\end{picture}\endgroup      \caption{$\mathrm{SPX}(r,1)) \cong \Inf[1][0][C_{2r}][][2]$}
  \end{subfigure}
  \caption{}\label{fig:PXandSPX}
\end{figure}

\begin{example}
  Natural examples of the graphs listed above are the Praeger-Xu graphs and the Split Praeger-Xu graphs (see \cite{PotocnikSpiga_2021_NumberFixedPoints}).
  Namely 
  $\mathrm{PX}(r,1) \cong \lex{(2 K_{1})}{C_{r}}$ and $\mathrm{SPX}(r,1)) \cong \Inf[1][0][C_{2r}][][2]$ with $\cP$ the partition taking every second edge (see \cref{fig:PXandSPX}).
\end{example}

In order to prove \cref{thm:graphs} we obtain a number of results regarding the minimal degree of permutation groups which might be of interest in their own right.
 The study of minimal degree of permutation groups goes all the way to Jordan \cite{Jordan_1871_TheoremesSurLes}
(see also \cite[Thm. 3.3D]{DixonMortimer_1996_PermutationGroups}),
who proved that for every constant $C$, every primitive permutation group with minimal degree at most $C$ is either symmetric, alternating or it belongs to a finite list of exceptions.
This result has been improved by several authors (see
\cite{ GuralnickMagaard_1998_MinimalDegreePrimitive, LiebeckShalev_2015_FixedPointsElements, BurnessGuralnick_2022_FixedPointRatios}
for example) and as a consequence, primitive permutation groups $G \leq \sym(\Omega)$ with $\mindeg (G)< \frac{2|\Omega|}{3}$ are completely classified.

Much less is known about the minimal degree of imprimitive permutation groups.
Even though some work has been done in this context (see \cite{SaxlShalev_1995_FixityPermutationGroups,LiebeckSaxl_1991_MinimalDegreesPrimitive,LawtherLiebeckSeitz_2002_FixedPointRatios,Burness_2007_FixedPointRatios1}, for example), several basic questions, such as \cref{prob:theproblem} below, remain unresolved.

\begin{problem}\label{prob:theproblem}
Given a (small) constant $d$, determine all transitive permutation groups of minimal degree $d$.
\end{problem}

In this paper we address this problem when $d$ is an arbitrary prime number (see \cref{thm:groups_minDegp})  as well as when $d=4$ (\cref{thm:groups_minDeg4}).
Our analysis of the case $d=4$  suggests that for composite integers $d$ the complexity of  \cref{prob:theproblem} grows very quickly with $d$.

In order to state \cref{thm:groups_minDegp} and \cref{thm:groups_minDeg4} we need to introduce some additional notation.
For a permutation group $G \leq \sym(\Delta)$ and a set $\Lambda$, let $\fun(\Lambda,G)$ be the group of all functions $\bfg :\Lambda \to G$ (with the product defined point-wise).
The group $\fun(\Lambda,G)$ then acts on the set $\Omega = \Delta \times \Lambda$ according to the rule \[  (\delta,\lambda)^{\bfg} = (\delta^{\bfg(\lambda)}, \lambda).\]
Suppose in addition that $H$ is a permutation group on $\Lambda$. Then every element $h \in H$ induces an automorphism of $\fun(\Lambda,G)$ given by \[\bfg^{h} (\lambda) = \bfg(\lambda^{h^{-1}}).\]

If $\Lambda$ is finite, say $|\Lambda|=k$, and its elements are labelled by $1, \dots, k$, we may identify the function $\bfg$ with a $k$-tuple $(g_{{1}}, \dots, g_{{k}}) \in G^{k}$, and thus the group $\fun(\Lambda, G)$ with  the direct product $G^{k}$ of $k$ copies of $G$. In this context, the action of $H$ on $G^{k}$ is given by permuting the coordinates of each of the $k$-tuples.

The action of $H$ on $\fun(G,\Lambda)$ defines the \emph{(imprimitive) wreath product}
\[G \wr H := \fun(\Lambda, G) \rtimes H\]
acting on $\Delta \times \Lambda$ by
\[(\delta,\lambda)^{(\bfg, h )} = \left( \delta^{\bfg(\lambda)}, \lambda^{h} \right).\]
The group $\fun(\Lambda,G)$ is called the \emph{base group} of $G \wr H$.

Following \cite{PraegerSchneider_2018_PermutationGroupsCartesian}, if $G_{1}$ is a permutation group on $\Omega_{1}$ and $G_{2}$ is a permutation group on $\Omega_{2}$, we say that a pair $(f,\phi)$ with $f:\Omega_{1} \to \Omega_{2}$ a bijection and $\phi: G_{1} \to G_{2}$ an injective homomorphism, is a \emph{permutation embedding} if
\begin{equation}\label{eq:permEmb}
  f(\omega^{x}) = \left( f(\omega) \right)^{ \phi (x)}
\end{equation}
for every $\omega \in \Omega_{1}$ and $x \in G_{1}$.
\cref{eq:permEmb} is equivalent to saying that the diagram
\[
  \begin{tikzcd}
    \Omega_{1} \arrow[d, "x"'] \arrow[r, "f"] & \Omega_2 \arrow[d, "\phi (x)"] \\
    \Omega_{1} \arrow[r, "f"]                 & \Omega_{2}
  \end{tikzcd}
  \]
  is commutative for every $x \in G_{1}$.
In addition, if $\phi$ is a group isomorphism then the pair $(f, \phi)$ is a \emph{permutation isomorphism} and the permutation groups $G_{1}$ and $G_{2}$ are \emph{permutation isomorphic}. In this paper, we often abuse notation and simply write $G_{1} \leq G_{2}$ and $G_{1} \cong G_{2}$ to indicate that $G_{1}$ is permutation isomorphic to a subgroup of $G_{2}$ and that $G_{1}$ is permutation isomorphic to $G_{2}$, respectively.
  
For $k \in \bN$, we denote by $[k]$ the set $\left\{ 1, \dots, k \right\}$. 
  We are now ready to state our main results regarding transitive permutation groups with small minimal degree. 

\begin{thm}\label{thm:groups_minDegp} 
  Let $p$ be a prime number and $G$ a transitive permutation group of degree $n$ containing a $p$-cycle.
  Then there exist $m,k \in \bN$ satisfying $n=mk$, $p \leq m \leq n$, and permutation groups $X \leq Y$  of degree $m$, such that:
  \begin{enumerate}
     \item \label{item:groups_minDegP_wr}  up to permutation isomorphism,  $G \leq Y \wr \sym(k)$, acting on $ [m]\times [k] $;
     \item \label{item:groups_minDegP_base} the intersection of $G$ with the base group $Y^{k}$ of $ Y \wr \sym(k)$ contains $X^{k}$;
    \item \label{item:groups_minDegP_table} $p,m,X$ and $Y$ are as in \cref{tab:degreeP}.
    \end{enumerate}
   Moreover, the minimal degree of $G$ is $p$ if and only if the condition in the last column of \cref{tab:degreeP} holds.

\end{thm}

\begin{table}[htb]
  \begin{tabularx}{\textwidth}{|\mm{1.8}|\mm{.6}|\mm{.8}|\mm{.8}|\mm{1}|}
    \hline
    p & m &X & Y & \mindeg(G) = p?\\ \hline
2 & \geq 2 & \sym(m) & \sym(m)& \text{always} \\
     \geq 3 & \geq 3 & \alt(m) & \sym(m) & p=3 \text{ and (C)} \\
     \geq 5   & p  & C_{p} & \agl_{1}(p) & \text{(C)} \\
     \frac{q^d - 1}{q-1} & p & \pgl_{d}(q) & \pGl_{d}(q) & \text{never} \\
     11 & 11 & \psl_{2}(11) & \psl_{2}(11) & \text{never}\\
     11 & 11 & M_{11} & M_{11} & \text{never} \\
     23 & 23 & M_{23} & M_{23} & \text{never} \\
     \text{Mersenne prime } 2^{a}-1 & 2^{a} &\agl_{1}(2^{a}) & \aGl_{1}(2^{a}) & \text{(C)}  \\
     \text{Mersenne prime } 2^{d}-1 & 2^{d}  &\agl_{d}(2) & \aGl_{d}(2) & \text{never}  \\
     \geq 5 & p+1 & \psl_{2}(p) & \pgl_{2}(p) & \text{never} \\
     11 & 12 & M_{11} & M_{11} & \text{never} \\
     11 & 12 & M_{12} & M_{12} & \text{never} \\
     23 & 24 & M_{24} & M_{24} & \text{never} \\
\text{Mersenne prime } 2^{a}-1 & 2^{a}+1 & \pgl_{2}(2^{a}) & \pGl_{2}(2^{a}) & \text{(C)}\\

      \hline
  \end{tabularx}
  \caption{Permutation groups containing a $p$-cycle. Condition (C) appearing in the last column means that the pointwise stabiliser $G_{(\Omega \sm ([m]\times \{1\}) )}$ of  all points with second coordinate different from  $1$ is permutation isomorphic to $X$.}
  \label{tab:degreeP}\end{table}

\begin{remark}
  Note that since the value $\frac{q^{d}-1}{q-1}$ in row $4$ of \cref{tab:degreeP} is a prime, one can show that $\pgl_{d}(q)=\psl_{d}(q)$, implying that
   $\pgl_{d}(q)$ is simple.
\end{remark}

\begin{remark}\label{rem:cfsg}
  The proof of \cref{thm:groups_minDegp} relies on the results developed by Jones in \cite{Jones_2014_PrimitivePermutationGroups}, which in turn depend on the classification of finite simple groups (CFSG). 
  However, \cref{item:graphs_mindeg2_,item:graphs_mindegp_} of \cref{thm:graphs}, the graph theoretical results derived from \cref{thm:groups_minDegp}, can be proved without the assumption of the CFSG (see \cref{sec:graphs}). 
\end{remark}

\begin{rem}
  \cref{thm:groups_minDegp} offers a complete classification of transitive permutation groups with minimal degree $2$. 
  Such a group has to be permutation isomorphic to $\sym(m) \wr H$ for a transitive group $H \leq \sym(k)$.
  However, even for transitive permutation groups $G$ of minimal degree $3$ all we can deduce from \cref{thm:groups_minDegp} is that 
  \begin{equation} \label{eq:minDeg3}
    \alt(m)^{k} \leq G\leq \sym(m) \wr \sym(k).  
  \end{equation}
  While this description of transitive permutation groups of minimal degree $3$ is somewhat unsatisfactory from a group theoretical viewpoint, it suffices to classify vertex-transitive graphs whose automorphism group contains a $3$-cycle (see \cref{thm:graphs}). This is why we will not pursue a more detailed classification of such groups but rather pose it as a challenge for future work:
\end{rem}

\begin{prob}
  Provide a detailed description of transitive permutation group of minimal degree $3$. In particular, determine which groups $G$ satisfying \eqref{eq:minDeg3} are transitive.
\end{prob}

Let us now move our attention to transitive groups of minimal degree $4$.
Note that every permutation moving only $4$ points is either a $4$-cycle or a product of two disjoint transpositions; 
we shall call the latter a \emph{$2^{2}$-element}. 
Since the square of a $4$-cycle is a $2^{2}$-element, every permutation group of minimal degree $4$ contains a $2^{2}$-element.

\begin{thm}\label{thm:groups_minDeg4}
  Let $G$ be a transitive permutation group of degree $n$ containing a $2^{2}$-element. Then, up to permutation isomorphism, $G$ satisfies one of the following:
  \begin{enumerate}
    \item  \label{item:groups_minDeg4_prim} There exist $m,k \in \bN $ satisfying $n = mk$, $4 \leq m \leq n$, and  permutation groups $X \leq Y$ of degree $m$ such that
      \begin{enumerate}
        \item \label{item:groups_minDeg4_prim_wr} $G \leq Y \wr \sym(k)$, acting on $[m] \times [k]$;
        \item \label{item:groups_minDeg4_prim_base} the intersection of $G$ with the base group $Y^{k} $ of $ Y \wr \sym(k)$ contains $X^{k}$; and
        \item \label{item:groups_minDeg4_prim_tab} $m$, $X$ and $Y$ are as in \cref{tab:22-prim}.
      \end{enumerate} 
    In this case, $\mindeg(G) = 4$ unless $(X,Y) = (\alt(m),\sym(m))$, when $\mindeg(G) \leq 3$.

    \item \label{item:groups_minDeg4_cross} There exist $m,k \in \bN$ satisfying $n=2mk$, and  permutation groups $X \leq Y \leq \sym(2) \wr \sym(m)$ acting on the set $\Delta = [2] \times [m]$ such that:
\begin{enumerate}
      \item \label{item:groups_minDeg4_cross_wr} $G \leq Y \wr \sym(k)$, acting on the set $\Delta \times [k]$;
      \item \label{item:groups_minDeg4_cross_base} the intersection of $G$ with the base group $Y^{k}$ of  $ Y \wr \sym(k)$ contains $X^{k}$; and
      \item \label{item:groups_minDeg4_cross_tab} $X$ and $Y$ are as in \cref{tab:22-cross}.
    \end{enumerate}

    \end{enumerate}
\end{thm}

\begin{table}[hbt]
  \begin{tabularx}{.8\textwidth}{|M|M|M|}
    \hline
    m & X & Y \\ \hline
    \geq 4 & \alt(m) & \sym(m)\\
     5 & \dih(5) & \agl_{1}(5) \\
     6 & \psl_{2}(5) & \pgl_{2}(5) \\
     7 & \pgl_{3}(2)  & \pgl_{3}(2) \\
     8 & \agl_{3}(2) & \agl_{3}(2) \\ \hline
  \end{tabularx}
  \caption{Primitive groups containing a $2^{2}$-element. }
  \label{tab:22-prim}
\end{table}

\begin{table}[hbt]
  \begin{tabularx}{.8\textwidth}{|M|M|}
    \hline
     X & Y \\ \hline
\left\{ 1 \right\} \times \sym(m) & \left\langle \sflip \right\rangle \times \sym(m)  \\
(\sym(2))^{+} \times \left\{ 1 \right\}  & \sym(2) \wr \sym(m) \\ \hline
\end{tabularx}
  \caption{Subgroups of $\sym(2) \wr \sym(m)$ containing a $2^{2}$-element. Here $\sflip = ( (1\ 2 ), \dots, (1\ 2)) \in \sym(2)^{m}$ is the diagonal of the base group of $\sym(2) \wr \sym(m)$ and $(\sym(2))^{+}$ is the ``augmentation'' subgroup of $\sym(2)^{m}$ consisting of all the $m$-tuples that are non-trivial in an even number of entries.}
  \label{tab:22-cross}
\end{table}

We prove \cref{thm:groups_minDegp} and \cref{thm:groups_minDeg4} in \cref{sec:groups} while \cref{thm:graphs} is proved in \cref{sec:graphs}.

 \section{Permutation groups with small minimal degree} \label{sec:groups}

Before going into the proofs of our results let us establish some notation and terminology.

For a permutation group $G \leq \sym(\Omega)$ and a set $B \subseteq \Omega$, we let $G_B = \{g\in G : B^g =B\} $ denote the set-wise stabiliser of $B$,
and let $G_{(B)} = \cap_{\omega\in B} G_\omega$ be the  point-wise stabiliser of $B$.
If $G$ acts upon a set $X$, then we let $G^X$ denote the permutation group induced by this action.
In particular, if $B \subset \Omega$, then $G_B^B$ denotes the permutation group induced by the action of $G_B$ on $B$.
Note also that the pointwise stabiliser $G_{(\Omega \sm B)}$ of the complement of $B$ in $\Omega$ acts faithfully on $B$ and can thus be thought of as a subgroup of $G_{B}^{B}$.  

If $G$ is a transitive permutation group on  $\Omega$
we define a \emph{block system} for $G$ to be
a $G$-invariant partition of $\Omega$ consisting of sets of cardinality at least $2$; the elements of the block system are then called \emph{blocks}.
If $\cB$ and $\cD$ are block systems for $G$ such that every block of $\cB$ is contained in a block of $\cD$
then we say that $\cB$ is \emph{finer} than $\cD$ (or that $\cD$ is \emph{coarser} than $\cB$).
A \emph{system of minimal blocks} is a block system $\cB$ such that the only block system finer than $\cB$ is $\cB$ itself.
Note that $G$ is \emph{primitive} if a system of minimal blocks for $G$ consists of one block only, that is, if the only block system of $G$ is $\{\Omega\}$.
Otherwise, $G$ is \emph{imprimitive.}
Observe that if $B$ is a minimal block, then $G_B^B$ is a primitive permutation group.

Wreath products and imprimitive groups are naturally connected (see \cite[Section 5.2]{PraegerSchneider_2018_PermutationGroupsCartesian} and \cite[Section 2.6]{DixonMortimer_1996_PermutationGroups}).
The following theorem is a classical result; a proof can be found in \cite[Thm. 5.5]{PraegerSchneider_2018_PermutationGroupsCartesian}, see also \cite[Theorem 1.8]{Cameron_1999_PermutationGroups}.

\begin{thm} \label{thm:emb_wr}
  Let $G$ be an imprimitive group on $\Omega$. 
  Let $\cB$ be a block system for $G$, and let $B \in \cB$.
  Then the group $G$ is permutation isomorphic to a subgroup of $G_{B}^{B} \wr G^{\cB}$ acting naturally on $B \times \cB$.
  Moreover,  the subgroup $\left\langle  G_{(\Omega \sm B)} : B \in \cB \right\rangle$ of $G$ is isomorphic to $ \prod_{B \in \cB} G_{(\Omega \sm B)}$ and embeds naturally into the base group of $G_{B}^{B} \wr G^{\cB}$.
\end{thm}

In light of \cref{thm:emb_wr}, if $G$ is a permutation group with a block system $\cB$ and if $B \in \cB$ and $X \leq G_{(\Omega \sm B)}$, then $G$ contains a group isomorphic to $X^{|\cB|}$ that is embedded into the base group of $G_{B}^{B} \wr G^{\cB}$.
We may thus abuse notation and simply write 
\begin{equation}\label{eq:lower}
  X^{|\cB|} \leq G.
\end{equation}
Conversely, whenever we see the expression \eqref{eq:lower}  we tacitly assume that $G$ is a permutation group with a block system $\cB$ and that $X \leq  G_{(\Omega\sm B)}$ for a block $B \in \cB$.

A particular instance of the situation above is stated in the following lemma which will be used frequently in the development of our results.

\begin{lemma}\label{lem:theOne}
  Let $G$ be a permutation group on $\Omega$ and let $\cB$ be a system of minimal blocks for $G$.
  Assume that $x \in G$ is such that $\supp(x) \subset B$ for some $B \in \cB$.
  Let $X = \left\langle x^{G_{B}} \right\rangle = \left\langle x^{g} : g\in G_{B} \right\rangle $. Then up to permutation isomorphism,  \[ X^{|\cB|} \leq G \leq G_{B}^{B} \wr G^{\cB}\] with $X^{|\cB|}$ contained in the base group of $G_{B}^{B} \wr G^{\cB}$.
\end{lemma}

\subsection{Transitive permutation groups of prime minimal degree}

This section is devoted to the proof of \cref{thm:groups_minDegp}. 

We leave it to the reader to check that a transitive group $G$ satisfying conditions (\ref{item:groups_minDegP_wr}), (\ref{item:groups_minDegP_base}), and (\ref{item:groups_minDegP_table}) in \cref{thm:groups_minDegp} has minimal degree $p$ if and only if the condition in the last column of \cref{tab:degreeP} is fulfilled.

We shall thus assume throughout this section that $G$ is a transitive permutation group  containing a cycle of prime length $p$ and prove that $G$ satisfies conditions (\ref{item:groups_minDegP_wr}), (\ref{item:groups_minDegP_base}), and (\ref{item:groups_minDegP_table}) in \cref{thm:groups_minDegp}.

If $G$ is primitive, then the validity of \cref{thm:groups_minDegp} (with $k=1$ and $n=m$) can be checked using \cref{thm:primitiveP} below, which follows directly from the work of Gareth Jones; see \cite[Theorem 1.1 and Theorem 1.2]{Jones_2014_PrimitivePermutationGroups}.

\begin{thm}[{\cite[Theorem 1.1 and Theorem 1.2]{Jones_2014_PrimitivePermutationGroups}}]
  \label{thm:primitiveP}
Let $G$ be a primitive group of finite degree $m$. 
Let $p$ be a prime and suppose that $G$ contains a $p$-cycle fixing $k$ points (so that $m=p+k$).
Then one of the following holds:
\begin{enumerate}
  \item \label{item:primitiveP_alt}
$\alt(m) \leq G \leq \sym(m)$; or
  \item \label{item:primitiveP_0} $k = 0$, $p \geq 5$ and either:
    \begin{enumerate}
      \item \label{item:primitiveP_0a} $C_{p} \leq G \leq \agl_{1}(p)$; or
      \item \label{item:primitiveP_0b}$\pgl_{d}(q) \leq G \leq \pGl_{d}(q)$ acting on the points of the $(d-1)$-dimensional projective space over $\bF_{q}$, that is $m = q^{d-1} + q^{d-2} + \cdots + q + 1$ for some $d \geq 2$ and $q$ a prime power; or
      \item \label{item:primitiveP_0c} $G=\psl_{2}(11)$ acting on a Payley biplane with $11$ points, $M_{11}$ or $M_{23}$ acting on $11$ and $23$ points, respectively.
    \end{enumerate}
  \item \label{item:primitiveP_1} $k = 1$, $p \geq 5$ and either:
    \begin{enumerate}
      \item \label{item:primitiveP_1a}$\agl_{d}(2^{a}) \leq G \leq \aGl_{d}(2^{a})$ acting on the vector space $\bF_{2^{a}}^{d}$ and $p = 2^{ad}-1$ is a Mersenne prime (in particular $a=1$ or $d=1$); or
      \item \label{item:primitiveP_1b}$G = \psl_{2}(p)$ or $G=\pgl_{2}(p)$ acting on the projective line; or
      \item \label{item:primitiveP_1c} $G = M_{11}$ acting on $12$ points, or $G$ is either $M_{12}$ or $M_{24}$  with its natural action on $12$ and $24$ points, respectively.
    \end{enumerate}
  \item \label{item:primitiveP_2} $k=2$, $p \geq 5$ and $\pgl_{2}(2^{a}) \leq G \leq \pGl_{2}(2^{a})$ acting on the projective line (that is $m = 2^{a}+1$) and $p = 2^{a}-1$ is a Mersenne prime.
\end{enumerate}
\end{thm}

Now assume that $G$ is an imprimitive transitive permutation group on $\Omega$ containing a $p$-cycle $x = (\alpha_{1} \cdots  \alpha_{p})$ for some prime $p$.
Let $\cB$ be a system of minimal blocks for $G$.
Let $B \in \cB$ be the block containing $\alpha_{1}$ and let $m=|B|$. 
Then  $\supp(x) \subset B$; in particular we can think of $x$ as a nontrivial element of $G_{B}^{B}$. 

Let $Y = G_{B}^{B}$ and let $X = \left\langle x^{Y} \right\rangle $. 
Since $Y$ is primitive and $x \in Y$, it follows that $Y$ is one of the groups listed in \cref{thm:primitiveP}.
Furthermore, since $X \trianglelefteq Y$ and $Y$ is primitive, we see that $X$ is transitive on $B$.
In view of \cref{lem:theOne}, in order to prove \cref{thm:groups_minDegp}, we need to show that the pair $(X,Y)$  is as in \cref{tab:degreeP}.

Suppose that $Y$ is as in \cref{item:primitiveP_alt} of \cref{thm:primitiveP}; that is, $\alt(m) \leq Y \leq \sym(m)$. 
Since $X \trianglelefteq Y$, then $\alt(m) \leq X$.
If $p=2$ then  $X$ contains an odd permutation, which implies that $X=\sym(m)$. If $p$ is odd, then $X$ is generated by cycles of odd length, which implies that $X = \alt(m)$.
Therefore, $(X,Y)$ is as in one of the first two rows of \cref{tab:degreeP}.

Now let $C_{p} \leq Y \leq \agl_{1}(p)$.
Recall that $\agl_{1}(p)=C_{p}\rtimes \bF_{p}^{*} \cong C_{p} \rtimes C_{p-1}$.
Sylow's Theorem implies that the $p$-cycle $x$ is a conjugate of an element in $C_{p}$, however $C_{p} \trianglelefteq G$.
It follows that $X = C_{p}$, as in the third row of \cref{tab:degreeP}.

Assume now that $\pgl_{d}(q) \leq Y \leq \pGl_{d}(q)$ acting on the $\frac{q^{d}-1}{q-1}$ points of the projective space, as in \cref{item:primitiveP_0b} of \cref{thm:primitiveP}. 
Note that $\frac{q^{d}-1}{q-1}$ is assumed to be a prime number $p$, implying that $d$ is a prime and that $\psl_{d}(q) = \pgl_{d}(q)$. 
Since $X \trianglelefteq Y$ and the socle of $Y$ is $\pgl_{d}(q)$, it follows that $X = \pgl_{d}(q)$, as claimed in the fourth row of \cref{tab:degreeP}.
Moreover, it can be shown that every $p$-cycle in $Y$ is a conjugate of a Singer cycle in $\pgl_{d}(q)$ (see \cite[Corollary 2]{Jones_2002_CyclicRegularSubgroups}).

Suppose that $\agl_{d}(2^{a}) \leq Y \leq \aGl_{d}(2^{a})$.
We claim that $X=\agl_{d}(2^{a})$.
The fact that $\agl_{d}(2^{a}) \leq X$ is a consequence of \cite[Lemma 6.3]{Mueller_1996_ReducibilityBehaviorPolynomials}.
It is easy to see that the largest power of $p$ that divides $|\aGl_{d}(2^{a})|$ is $p$ itself hence any two $p$-cycles are conjugate.
The group $\agl_{d}(2^{a})$ is normal in $\aGl_{d}(2^{a})$ (and hence normal in $Y$) and contains a $p$-cycle (a Singer cycle in $\gl_{d}(2^{a})$) hence every $p$-cycle is contained in $\agl_{d}(2^{a})$.

A similar argument to the one above proves that if $Y$ is a group as in \cref{item:primitiveP_1b} or \cref{item:primitiveP_2}, then $X = \psl_{2}(q)$.
If $q = 2^{a}$ as in \cref{item:primitiveP_2} of \cref{thm:primitiveP}, then $\psl_{2}(q) = \pgl_{2}(q)$.

For the remaining cases, that is when $Y$ is as in \cref{item:primitiveP_0c,item:primitiveP_1c} of \cref{thm:primitiveP}, the group $Y$ is simple which implies that $X = Y$.

The discussion above together with \cref{lem:theOne} prove that if $G$ is a group containing a $p$-cycle, then $G$ is permutationally isomorphic to a subgroup of $Y \wr G^{\cB}$ whose intersection with the base group $Y^{|\cB|}$ contains $X^{|\cB|}$.
This of course is equivalent to the statement of \cref{thm:groups_minDegp}.

\subsection{Permutation groups with minimal degree 4}\label{sec:grpsMinDeg4}

The corresponding analysis for a permutation group containing an element fixing all but four points is slightly more complicated.
As mentioned before, it is enough for us to understand transitive permutation groups containing a $2^{2}$-element.
Let us first observe a basic structural result.
\begin{lem}\label{rem:deg4_supp}
  Let $G$ be a transitive permutation group on $\Omega$ containing a $2^{2}$-element $x=(\alpha_{1}\ \alpha_{2})(\beta_{1} \ \beta_{2})$. 
  Let $\cB$ be a system of minimal blocks for $G$ and let $B_{1}$ be the block containing $\alpha_{1}$.
  Then, up to swapping the roles of $\beta_{1}$ and $\beta_{2}$, one of the following holds:
  \begin{enumerate}
    \item \label{item:deg4_supp4} $\supp(x) \subset B_{1}$; or
    \item \label{item:deg4_supp2e} $B_{1} = \left\{ \alpha_{1}, \beta_{1} \right\} $.
    In this case all the blocks of $\cB$ are of size $2$ and the block $B_{2}$ containing $\alpha_{2}$ is precisely $\left\{ \alpha_{2},\beta_{2} \right\} $; or
    \item \label{item:deg4_supp2h} $\supp(x) \cap B_{1} = \left\{ \alpha_{1}, \alpha_{2} \right\} $, in which case there is another block $B_{2} \neq B_{1}$ that contains $\left\{ \beta_{1}, \beta_{2} \right\} $.
  \end{enumerate}
\end{lem}

We will now prove \cref{thm:groups_minDeg4} by considering each of the three cases mentioned above separately.
We shall first deal with the case where $\supp(x) \subset B_{1}$.
The analysis in this case is very similar to the one done for prime minimal degree.
If $\cB$ is a system of minimal blocks, then for each block $B \in \cB$, the action of $G_{B}$ on $B$ is primitive and the group $X=\left\langle x^{G_{B}} \right\rangle $ is contained in $G_{(\Omega \sm B)}$.
In this case, \cref{lem:theOne} implies that $G$ satisfies \[X^{|\cB|} \pleq G \pleq G_{B}^{B} \wr G^{\cB}. \]

We thus need to find all  possible pairs $(X,Y)$ with $Y$ a primitive permutation group cointaining a $2^{2}$-element $x$ and $X = \left\langle x^{Y} \right\rangle  \leq Y$.
The following classical results impose strong conditions on the group $Y$:

\begin{lem}[{\cite[Thm. 3.3D]{DixonMortimer_1996_PermutationGroups}}]\label{lem:DTSmallDeg}
  Let $G$ be a primitive permutation group on $\Omega$ containing an  element fixing all but $4$ points.
  If $G$ is $2$-transitive and $\left| \Omega \right| > 10$, then $G$ contains  $\alt(\Omega)$.
\end{lem}

\begin{lem}[see {\cite[Thm. 0.3]{Babai_1981_OrderUniprimitivePermutation}}]\label{lem:babai}
   Let $G$ be a primitive permutation group on $\Omega$ containing a  $2^{2}$-element.
   If $\left| \Omega \right| > 49 $, then $G$ is $2$-transitive.
\end{lem}

It follows that every primitive permutation group $Y$ of degree larger than $49$ containing a $2^{2}$-element $x$ must be either alternating or symmetric.
In this case, $X = \left\langle x^{Y} \right\rangle = \alt(\Omega)$.
On the other hand, if $|Y| \leq 48$, the possibilities for $X$ for each such $Y$ can be easily determined using GAP \cite{TheGAPGroup_2022_GapGroupsAlgorithmsGapGroupsAlgorithms} and are precisely those in \cref{tab:22-prim} (see also \cite[Example 3.3.1]{DixonMortimer_1996_PermutationGroups}).
This discussion proves that the claim in \cref{item:groups_minDeg4_prim}  of \cref{thm:groups_minDeg4} holds whenever $\supp(x) \subset B_{1}$.

\begin{figure}
  \def\svgwidth{.6\textwidth}
  \begin{scriptsize}
    \begingroup \makeatletter \providecommand\color[2][]{\errmessage{(Inkscape) Color is used for the text in Inkscape, but the package 'color.sty' is not loaded}\renewcommand\color[2][]{}}\providecommand\transparent[1]{\errmessage{(Inkscape) Transparency is used (non-zero) for the text in Inkscape, but the package 'transparent.sty' is not loaded}\renewcommand\transparent[1]{}}\providecommand\rotatebox[2]{#2}\newcommand*\fsize{\dimexpr\f@size pt\relax}\newcommand*\lineheight[1]{\fontsize{\fsize}{#1\fsize}\selectfont}\ifx\svgwidth\undefined \setlength{\unitlength}{126.14173228bp}\ifx\svgscale\undefined \relax \else \setlength{\unitlength}{\unitlength * \real{\svgscale}}\fi \else \setlength{\unitlength}{\svgwidth}\fi \global\let\svgwidth\undefined \global\let\svgscale\undefined \makeatother \begin{picture}(1,0.66349402)\lineheight{1}\setlength\tabcolsep{0pt}\put(0,0){\includegraphics[width=\unitlength,page=1]{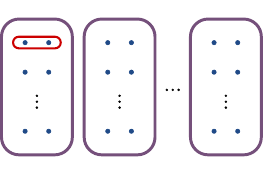}}\put(0.14044943,0.63604164){\color[rgb]{0.1254902,0.29019608,0.52941176}\makebox(0,0)[t]{\lineheight{1.25}\smash{\begin{tabular}[t]{c}$D_1$\end{tabular}}}}\put(0.14044943,0.00682824){\color[rgb]{0.8,0,0}\makebox(0,0)[t]{\lineheight{1.25}\smash{\begin{tabular}[t]{c}$\bar{D}_1$\end{tabular}}}}\put(0.45505619,0.63604164){\color[rgb]{0.1254902,0.29019608,0.52941176}\makebox(0,0)[t]{\lineheight{1.25}\smash{\begin{tabular}[t]{c}$D_2$\end{tabular}}}}\put(0.45505619,0.00682824){\color[rgb]{0.8,0,0}\makebox(0,0)[t]{\lineheight{1.25}\smash{\begin{tabular}[t]{c}$\bar{D}_2$\end{tabular}}}}\put(0.85955057,0.63604164){\color[rgb]{0.1254902,0.29019608,0.52941176}\makebox(0,0)[t]{\lineheight{1.25}\smash{\begin{tabular}[t]{c}$D_k$\end{tabular}}}}\put(0.85955057,0.00682824){\color[rgb]{0.8,0,0}\makebox(0,0)[t]{\lineheight{1.25}\smash{\begin{tabular}[t]{c}$\bar{D}_k$\end{tabular}}}}\put(0,0){\includegraphics[width=\unitlength,page=2]{blocks4easy.pdf}}\put(0.14044943,0.546154){\color[rgb]{0.1254902,0.29019608,0.52941176}\makebox(0,0)[t]{\lineheight{1.25}\smash{\begin{tabular}[t]{c}$B_1$\end{tabular}}}}\put(0.14044943,0.43379445){\color[rgb]{0.1254902,0.29019608,0.52941176}\makebox(0,0)[t]{\lineheight{1.25}\smash{\begin{tabular}[t]{c}$B_2$\end{tabular}}}}\put(0.14044943,0.20907543){\color[rgb]{0.1254902,0.29019608,0.52941176}\makebox(0,0)[t]{\lineheight{1.25}\smash{\begin{tabular}[t]{c}$B_m$\end{tabular}}}}\end{picture}\endgroup    \end{scriptsize}
  \caption{}
  \label{fig:blocks4easy}
\end{figure}

Assume now that $B_{1} = \left\{ \alpha_{1}, \beta_{1} \right\} $ as in \cref{item:deg4_supp2e} of \cref{rem:deg4_supp}.
Then every block in $\cB$ is of size $2$.
Let $\bar{\cD}=\left\{ \bar{D}_{1}, \dots, \bar{D}_{k} \right\}$ be a system of minimal blocks for the group $G^{\cB}$.
Observe that $\bar{\cD}$ induces a block system $\cD$ for $G$ that is coarser than $\cB$, where the blocks of $\cD$ are given by $D_{i} = \cup \bar{D}_i$. 
Observe that $B_{1}$ and $B_{2}$ belong to the same block of $\bar{\cD}$, which we may assume to be $D_{1}$ (see \cref{fig:blocks4easy}).

Let $D  = D_{1} \in \cD$, $\bar{D}= \bar{D}_{1}$ and $Y= G_{D}^{D}$.
Observe that $B_{1}=\left\{ \alpha_{1}, \beta_{1} \right\}, B_{2} = \left\{ \alpha_{2}, \beta_{2} \right\} \in \bar{D} $ and let $B_{3}, \dots, B_{m} \in \cB$ be such that $\bar{D}= \left\{ B_{1}, \dots, B_{m} \right\}$. 
\cref{lem:theOne} implies that, up to permutation isomorphism, $G\leq Y \wr \sym(k)$ 
and the group $Y$ is itself (permutation isomorphic to) a subgroup of $\sym(B_{1}) \wr \sym(\bar{D})$
acting on $B_{1} \times \bar{D}$.
To simplify notation we shall identify each pair $(\alpha_{1}, B_{i}) \in B_{1} \times \bar{D}$ with $i$ and the pair $(\beta_{1}, B_{i})$ with $-i$, 
and think of $Y$ as a permutation group on $\left\{ 1, -1, 2, -2, \dots, m,-m \right\}$.
Furthermore we shall identify $\sym(B_{1})$ with the multiplicative group $C_{2} = \left\{ 1, -1 \right\}$ and $\sym(\bar{D})$ with $\sym(m)$ so that $Y \leq C_{2} \wr \sym(m)$.
With these identifications, the $i^{th}$-coordinate of an element in $C_{2}^{m}$ acts upon $\left\{ i, -i \right\}$ by multiplication and the elements of $\sym(m)$ permute the blocks $B_i = \left\{ i,-i \right\}$, $i \in \left\{ 1, \dots, m \right\}$. \
Similarly, we shall identify $D$ with the set $\left\{ 1, -1, 2, -2, \dots, m,-m \right\}$.
Recall that $x$ swaps the blocks $B_{1}$ and $B_{2}$, and thus under the above identifications  induces the transposition $({1} \ {2}) \in Y^{\bar{D}}$. 
On the other hand, by the choice of $\bar{\cD}$,  $Y^{\bar{D}}$ is primitive, implying that $Y^{\bar{D}} = \sym(m)$.
Now, let $X = \left\langle x^{Y} \right\rangle$ and observe that $X \leq G_{(\Omega \sm D)}$. 
We shall now classify all the possibilities for $X$ and $Y$, which eventually will show that $G$ is as in \cref{item:groups_minDeg4_cross} of \cref{thm:groups_minDeg4}.
More preciselly we shall prove the following result.

\begin{prop}\label{prop:Easy}
  Let $G$ be a transitive permutation group acting on $\Omega$ with $|\Omega| = n$, admitting a block system $\cB$ with blocks of size $2$.
  Assume that $G$ contains a $2^{2}$-element $x = (\alpha_{1} \ \alpha_{2})(\beta_{1} \  \beta_{2})$ such that $\left\{ \alpha_{1}\ \beta_{1} \right\}$ and $\left\{ \alpha_{2}\ \beta_{2} \right\}$ are blocks in $\cB$.
  Then there exist $m,k \in \bN$, and permutation groups $X \leq Y \leq \cubeWr$, such that $n=2mk$, $(X,Y)$ are as in \cref{tab:22-easy}, and up to permutation isomorphism,
  \[G \leq Y \wr \sym(k),\]
  with intersection of $G$ with the base group of $Y \wr \sym(k)$ containing the group $X^{k}$.
  
  \begin{table}[hbt]
    \begin{tabularx}{.8\textwidth}{|M|M|}
      \hline
       X & Y \\ \hline
\left\{ 1 \right\} \times \sym(m) & \left\langle \sflip \right\rangle \times \sym(m)  \\
\evenAb \rtimes \sym(m) & \cubeWr \\
\hline
    \end{tabularx}
    \caption{Some transitive subgroups of $\cubeWr$ containing a $2^{2}$-element.}
    \label{tab:22-easy}
  \end{table}
\end{prop}

\begin{proof}
  As discussed above, the result follows from \cref{lem:theOne}.
  We just need to find the possibilities for $Y=G_{D}^{D}$ and $X = \left\langle x^{Y} \right\rangle $ as subgroups of $\cubeWr$ acting on $\left\{ \pm 1, \cdots, \pm m \right\}$ as discussed above.
  
  Consider the natural homomorphism $\pi: C_{2} \wr \sym(m) \to \sym(m)$ and notice that $\pi(Y) =Y^{\bar{D}} = \sym(m)$.
  Recall that $X$ is a normal subgroup of $Y$, which implies that $\pi(X)$ is a normal subgroup of $\pi(Y)$. 
  But since $\pi(X)$ contains a transposition (the one induced by $x$), it follows that $\pi(X) = \sym(m)$.
  Therefore, we need to classify the pairs $(X,Y)$ satisfying
  \begin{itemize}
    \item  $X \leq Y \leq \cube[m] $;
    \item  $X = \left\langle x^{Y} \right\rangle $  with $x$ a $2^{2}$-element such that $\pi(x)$ is a transposition; and
    \item $  \pi(X) =  \pi(Y) = \sym(m)$.
  \end{itemize}
  We shall classify them by considering the group $K = Y \cap \ker \pi \leq \ab[m]$.
  
  First, observe that $K$ is a normal subgroup of $Y$ contained in $\ab[m]$.
  It follows that $K$ must be preserved by $\pi(Y)=\sym(m)$.
  Those groups are easy to classify (see for example \cite[Lemma 3.1]{CollinsMontero_2021_EquivelarToroidsFew} or observe that $K$ can be viewed as a submodule of the permutation module of $\sym(m)$ over the field $\mathbb{F}_{2}$). In particular, one of the following must hold:
  \begin{align}\label{eq:kernel}
    K&=\left\{ \vid \right\},   &
    K&= \left\langle \sflip \right\rangle,  &
    K&= \evenAb[m], \ \text{or} &
    K&=  \ab[m],
  \end{align}
  where $\vid$ is the trivial element of $C_{2}^{m}$, $\sflip$ is the so called \emph{superflip}, that is, the permutation that swaps the two elements in each pair $\left\{ i, -i \right\}$; and
\[\evenAb[m] = \left\{ v \in \ab[m] : |\supp(v)| \text{ is even } \right\}, \]
  where $\supp(v)$ denotes the subset of $\{1, \dots, m\}$ where $v$ has non-trivial entry (often called the augmentation module).
  
  Before classifying the possibilities for $Y$ let us prove the following claims, which will impose certain restrictions for $X$:
  
\begin{enumerate}
      \item\label{item:ConjClass22_tr} If $X$ is transitive on $D = \left\{ 1, -1, 2, -2, \dots, m,-m \right\}$, then $\evenAb[m] \leq X$.
      \item\label{item:ConjClass22_abelPlus} If $m \geq 3$ and $\evenAb[m] \leq Y$, then $\evenAb[m] \leq X$.
      \item\label{item:ConjClass22_abel} If $\ab \leq Y$, then $\evenAb[m] \leq X$.
    \end{enumerate}
To prove Claim (\ref{item:ConjClass22_tr}), let us assume that $X$ is transitive on $D$.
      Since $\pi(Y) = \sym(m)$, for every $3 \leq i \leq m$ there exists $g_{i} \in Y$  that fixes $B_{1}$ and swaps $B_{i}$ and $B_{2}$. 
      Define $a_{1} =1$, $b_{1}= -1$ and $g_{2} = 1 \in Y$.
      Further, for $2 \leq i \leq m$, let $x_{i}:=x^{g_{i}}=(1 \ a_{i})(-1 \ b_{i})$ where $\left\{ a_{i}, b_{i} \right\} = \left\{ i, -i \right\}$.

      Let $X' = \left\langle x_{i} : 2\leq i \leq m \right\rangle  \leq X$.
      Observe that $X'$ has precisely two orbits on $D$, namely $\left\{ a_{i} : 1 \leq i \leq m \right\}$ and $\left\{ b_{i} : 1 \leq i \leq m \right\}$, and $X'$ acts as the symmetric group on each of these orbits.
      Since $X$ is transitive, there must exist $y \in X$, a conjugate of $x$, such that $y$ maps $a_{i}$ to $b_{j}$ (for some $1 \leq i < j \leq m$), and since $y$ is a $2^{2}$-element, then $y$ must be of the form $(a_{i}\ b_{j})(b_{i}\ a_{j})$.
      Observe that $x_{j}^{x_{i}} = (a_{i} \ a_{j})(b_{i} \ b_{j}) \in X'$ and $(x_{j}^{x_{i}})y = (a_{i}\ b_{i})(a_{j} \ b_{j}) \in \evenAb[m] \cap X$.
      Finally observe that the conjugates of $(a_{i}\ b_{i})(a_{j} \ b_{j})$ under $Y$  generate $\evenAb[m]$ (recall that $\pi(Y) = \sym(m)$).
      Since $X \trianglelefteq Y$ it follows that $\evenAb[m]\leq X$, which proves Claim (\ref{item:ConjClass22_tr}).
      
      Now assume that $m \geq 3$ and that $\evenAb[m] \leq Y$, or that $C_{2}^{m} \leq Y$.
      Then there exists  $g \in Y$ that  swaps $1$ and $-1$ and fixes both $2$ and $-2$.
      Note that $x x^{g}  $ is the permutation $ (1 \ -1)(2\ -2) \in X \cap \evenAb[m]$.
      As before, this implies that $\evenAb[m] \leq X$; which proves Claims (\ref{item:ConjClass22_abelPlus}) and (\ref{item:ConjClass22_abel}).

      Now we classify the pairs $(X,Y)$ depending on the four possibilities for $K=Y\cap \ab[m]$ in \eqref{eq:kernel}.
      
      Assume that $K$ is trivial. 
      In this case $\pi$ restricted to $Y$ is an isomorphism between $Y$ and $\sym(m)$.
      Since $\pi(X) = \sym(m)$ we see that $X=Y$. 
      In particular, $Y$ is a permutation group isomorphic to $\sym(m)$ acting transitively on $2m$ points and contains a $2^{2}$-element.
      Notice that this implies that $m \geq 3$.
It is not hard to see that $Y$ must be permutation isomorphic (under $\pi$) to $\sym(m)$ acting on the cosets of $\alt(m-1)$.
      This action admits a block system with two blocks, namely:
      \[\begin{aligned}
        \left\{ \alt(m-1) \sigma : \sigma \in \alt(m) \right\} && \text{and} &&
        \left\{ \alt(m-1) \sigma : \sigma \in \sym(m) \sm \alt(m) \right\}
      \end{aligned}\]
      with odd permutations in $\sym(m)$ swapping these blocks.
      
      Now recall that
      $\pi(x)$ is a transposition, and thus $\pi(x) \in \sym(m)\sm \alt(m)$, implying that it swaps the above two blocks.
      But then, $x$ has no fixed points in its action on the cosets of $\alt(m-1)$
      One the other hand, $x$ should fix all but four points, implying that $m=2$, which is a contradiction.
      This contraduction excludes the case $K=\{ \vid \}$.
      
      Assume now that $K= \left\langle \sflip \right\rangle$, with $\sflip$ the superflip.
Recall that $\pi(X)=\sym(m)$ hence $X$ has at most two orbits.
      However, Claim (\ref{item:ConjClass22_tr}) implies that if $X$ is transitive, then $ \left\langle \sflip \right\rangle = K =  \evenAb[m] $,  which in turn implies that $m=2$. 
      In this situation,
      it is not hard to see that $Y= K \times \sym(2)$. 
      If $m \geq 3$, then $X$ must have exactly two orbits and each block of $B_i$ contains precisely one element of each orbit. 
Up to a relabelling of the elements of each pair $B_{i}$ (and hence a permutation isomorphism), we may actually assume that the orbits of $X$ are $\left\{ 1, \dots, m \right\}$ and $\left\{ -1, \dots, -m \right\}$.
Since $X \trianglelefteq Y$, then its orbits are blocks of imprimitivity for $Y$, which are swapped by the superflip $\sflip$.
      It follows that $Y=K \times X \cong \left\langle \sflip \right\rangle \times \sym(m)$.

        If $K = \evenAb$, then $|Y|=|\ker\pi|| \pi(Y)| = 2^{m-1}m!$.
        Hence, $Y$ is an index-$2$ subgroup of $\cube$.
        Those groups are classified in \cite[Cor. 3.5]{CollinsMontero_2021_EquivelarToroidsFew}.
        In this situation $Y$ must be either $\evenAb[m] \rtimes \sym(m)$ or
\[\evencubeWr = \left\{ \left(v, \sigma\right) \in \cube[m] : (-1)^{|\supp(v)|} = \sgn(\sigma) \right\}. \]
        
        However, the latter is impossible because either $x= (1\ 2)(-1\ -2)$ or $x = (1\ -2)(-1\ 2)$, and none of them belong to $\evencubeWr$.
It follows that if $K = \evenAb$, then $Y= \evenAb \rtimes \sym(m)$.
        If $m \geq 3$ then Claim (\ref{item:ConjClass22_abelPlus}) implies that $\evenAb \leq X$, and since $\pi(X) = \sym(m)$, then $X=Y$.
        If $m=2$ just observe that $\evenAb[2]$ is the group generated by the super flip, $X = \left\langle x \right\rangle  = \sym(2) $, and this case reduces to a previous one.
        Finally if $K=\ab[m]$, $|Y|=|\ker\pi|| \pi(Y)| = 2^{m}m!$, which implies that $Y = \cube[m]$.
        Claim (\ref{item:ConjClass22_abel}) implies that $\evenAb \rtimes \sym(m) \leq X$ however $X$ cannot be $\cubeWr$ because $X$ consists of only even permutations.
      \end{proof}

Let us move to the case when $G$ and $\cB$ are as in \cref{item:deg4_supp2h} of \cref{rem:deg4_supp}.
More precisely, we shall prove the following result.

\begin{prop}\label{prop:Hard}
  Let $G$ be a transitive permutation group containing a  $2^{2}$-element $x = (\alpha_{1}\ \alpha_{2})(\beta_{1}\ \beta_{2})$. 
  Assume that $G$ admits a system of minimal blocks $\cB$ such that for some $B_{1}, B_{2} \in \cB$, $\supp(x) \cap B_1 = \left\{ \alpha_{2}, \alpha_{2} \right\}$ and $\supp(x) \cap B_{2} = \left\{ \beta_{1}, \beta_{2} \right\}$.
  Then one of the following hold:
  \begin{enumerate}
   \item \label{item:lessThan4} $\mu(G) < 4$;
   \item \label{item:4prim} $\mu(G) = 4$ and there is a $2^{2}$-element $x'$ such that $\supp(x') \subset B$ for some $B \in \cB$. In this case, $G$ is as in \cref{item:groups_minDeg4_prim} of \cref{thm:groups_minDeg4}; or
   \item \label{item:4easy} $\mu(G)=4$ and there is a block system $\cP$ for $G$ consisting of blocks of size $2$ such that $\left\{ \alpha_{1}, \beta_{1} \right\}$ and $\left\{ \alpha_{2}, \beta_{2} \right\}$ are blocks of $\cP$. 
   In this case, $G$ is as in \cref{prop:Easy};
   \item \label{item:4hard} $\mu(G) = 4$,  the blocks of $\cB$ are of size $2$ and there exists a block system $\cD$, coarser than $\cB$, satisfying $|\cD|=k$, for $D \in \cD$, $|D| = 2m$ (hence $|\Omega|=2mk$) and the group $Y=G_{D}^{D}$ is permutation isomorphic to a subgroup of $\cubeWr$ containing the subgroup $\evenAb$.
   In this case $G$ satisfies \[ G \leq Y \wr \sym(k) \] and the intersection of $G$ with the base group contains a subgroup isomorphic to $(\evenAb)^{k}$. 
  \end{enumerate}
\end{prop}

  \begin{proof}
    The strategy for the proof goes as follows: we shall assume that $G$ does not contain a transposition, a $3$-cycle, or a $2^{2}$-element whose support is contained in a block of $\cB$; that is, that none of \cref{item:lessThan4,item:4prim} hold.
    We will thus see that under these assumptions, if $B \in \cB$ satisfies that $|B| \geq 3$, then we can build a system of minimal blocks $\cP$ for $G$ with blocks of size $2$ such that $G$ satisfies \cref{item:4easy}. 
    Finally, we shall prove that the remaining cases satisfy \cref{item:4hard}. 
    We do this by prooving a series of mostly straighforward claims.
    
    Let $B_{1},B_{2} \in \cB$ such that $\left\{ \alpha_{1}, \alpha_{2} \right\} \subset B_{1}$ and $\left\{ \beta_{1}, \beta_{2} \right\} \subset B_{2}$. 
    Observe that since $\cB$ consists of minimal blocks, for every block $B$, the group $G_{B}^{B}$ is primitive and  contains a transposition (for example $(\alpha_{1} \ \alpha_{2}) \in G_{B_1}^{B_1}$). 
    It follows that $G_{B}^{B} = \sym(B)$.
    
    \begin{claim}\label{claim:exists}
      For every $\alpha \in B_{1}\sm \left\{ \alpha_{1}, \alpha_{2} \right\} $ there exist $\beta \in B_{2}\sm \left\{ \beta_{1}, \beta_{2} \right\} $ and $\beta_{0} \in \{\beta_{1},\beta_{2}\}$ such that $(\alpha_{1}\ \alpha) (\beta_{0}\ \beta) \in G$.
    \end{claim}
    
    \begin{proof}\renewcommand{\qedsymbol}{}
      Since $G_{B_{1}}^{B_{1}} = \sym(B)$, there exists $g \in G_{B}$ such that $\alpha_{1}^{g} = \alpha_{1}$ and $\alpha_{2}^{g} = \alpha$. Then $x^{g} = (\alpha_{1}^{g}\ \alpha_{2}^{g})(\beta_{1}^{g}\  \beta_{2}^{g})=(\alpha_{1}\ \alpha)(\beta_{1}^{g}\ \beta_{2}^{g} )$.
      If $|\{\beta_{1}^{g}, \beta_{2}^{g}\} \cap \left\{\beta_{1}, \beta_{2}  \right\}| \neq 1 $ then $(x x^{g})^{2}$ is a $3$-cycle with support $\left\{ \alpha_{1}, \alpha_{2}, \alpha \right\} $. 
      Take $\beta_{0} \in\{\beta_{1}^{g}, \beta_{2}^{g}\} \cap \left\{\beta_{1}, \beta_{2}  \right\}$.
      Finally, just observe that $\beta \in \left\{ \beta_{1}^{g}, \beta_{2}^{g} \right\} \sm \left\{ \beta_{0} \right\} $ must be an element of $B_{2}$.
    \end{proof}
    
    \begin{claim} \label{claim:pivot}
      If $\beta_{0} \in \left\{ \beta_{1}, \beta_{2} \right\} $ is such that for some $\alpha \in B_{1} \sm \{\alpha_{1},\alpha_{2}\}$ and $\beta \in B_{2} \sm \left\{ \beta_{1}, \beta_{2} \right\} $ the permutation $(\alpha_{1}\ \alpha)(\beta_{0}\ \beta ) \in G$, then $\beta_{0}$ is the only element in $\left\{ \beta_{1}, \beta_{2} \right\} $ with this property.
    \end{claim}
\begin{proof}\renewcommand{\qedsymbol}{} 
      Assume otherwise, that is, that there exist $\alpha, \alpha' \in B_{1}$ and $\beta, \beta'\in B_{2}$ such that $y_{1}=(\alpha_{1}\ \alpha)(\beta_{1}\ \beta)$ and $y_{2}= (\alpha_{1}\ \alpha')(\beta_{2}\ \beta')$ are both in $G$.
      If $\alpha = \alpha'$, then $y_{1}y_{2}$ is either a $3$-cycle or a $2^{2}$-element (depending on whether or not $\beta = \beta'$) with support contained in $B_{2}$.
      If $\alpha \neq \alpha'$ and $\beta \neq \beta'$, then $(y_{1}y_{2})^{2} = (\alpha'\ \alpha\ \alpha_{1})$ and if $\alpha \neq \alpha'$ but $\beta = \beta'$, then $(xy_{1}y_{2})^{2} = (\alpha_{1}\ \alpha) (\alpha_{2}\ \alpha')$.
      In any case, we get a contradiction to our assumptions on $G$. 
    \end{proof}
    
    The previous claim allows us to assume, without loss of generality, that $\beta_{0} = \beta_{1}$.
    
    \begin{claim} \label{claim:unique}
      For each $\alpha \in B_{1}$ there exists a unique $\beta \in B_{2}$ such that $(\alpha_{1}\ \alpha) (\beta_{1} \ \beta) \in G$.
    \end{claim}
    \begin{proof}\renewcommand{\qedsymbol}{}
      Assume otherwise.
      If $\alpha = \alpha_{1}$, then  $(\alpha_{1}\ \alpha)(\beta_{1}\ \beta)=(\beta_{1}\ \beta)$, and thus $\beta$ must be $\beta_{1}$.
      If $\alpha \neq \alpha_{1}$, assume that $\beta$ and $\beta'$ are such that $y_{1}=(\alpha_{1}\ \alpha) (\beta_{1} \ \beta)$ and $y_{2}=(\alpha_{1}\ \alpha) (\beta_{1} \ \beta' ) $. 
      Then $y_{1}y_{2} = (\beta_{1}\ \beta\ \beta')$. 
      Either way, this implies that $G$ satisfies \cref{item:lessThan4}.
    \end{proof}

    The previous analysis defines an equivalence relation on the blocks of $\cB$ given by $B \sim B'$ if $B=B'$ or if $B \neq B'$ and there is a $2^{2}$-element $y= (\gamma_{1}\ \gamma_{2})(\delta_{1}\ \delta_{2})$ with $\gamma_{1},\gamma_{2} \in B$ and $\delta_{1}, \delta_{2} \in B' $.
    In this case, $G$ has a subgroup acting simultaneously as $\sym(B)$ on $B$ and $B'$.
    More precisely, there is a bijection $\phi: B \to B'$ such that for every permutation $g \in \sym(B)$ the permutation $g^{\phi} = \phi^{-1} g \phi: B' \to B' $ satisfies that $g g^{\phi} \in G$.
    This relation is preserved under $G$, hence its equivalence classes induce a block system $\bar{\cD} = \left\{ \bar{D}_{1}, \dots, \bar{D}_{k} \right\}$ for $G^{\cB}$.
    The system $\bar{\cD}$ in turn induces a system $\cD = \left\{ D_{1}, \dots, D_{k} \right\}$ for $G$ by $D_{i} = \cup \bar{D}_{i}$.
    The system $\cD$ is coarser than $\cB$.
    Assume that $D_{1}$ is the block in $\cD$ that contains $\alpha_{1}$.

    \begin{claim}\label{claim:Blarge}
      If $  \ell := |B_{1}| \geq 3$, then $ D_{1} = B_{1} \cup B_{2}$.
    \end{claim}
    \begin{proof}\renewcommand{\qedsymbol}{}
      Indeed, if $\alpha_{1}, \alpha_{2}, \alpha_{3} \in B_{1}$ and $B_{3}$ is a block in $\cB$ other than $B_{1}$ and $B_{2}$, then there exists $\gamma_{1}, \gamma_{3} \in B_{3}$ such that $y=(\alpha_{1}\ \alpha_{3})(\gamma_{1}\ \gamma_{3}) \in G$ and since $x = (\alpha_{1}\ \alpha_{2})(\beta_{1}\ \beta_{2})$  then $(xy)^{2} = (\alpha_{1}\ \alpha_{3}\ \alpha_{2}) \in G$, which is a contradiction.
    \end{proof}
    
    Observe that in the case discussed above the group $\sym(\ell)$ acts simultaneously on the two blocks of $\cB$ contained in $D$ and it is actually a subgroup of $G_{(\Omega \sm D)}$.
    Label the elements of $D_{1}$ such that $B_{1} = \left\{ \alpha_{i} : 1 \leq i \leq \ell \right\} $, $B_{2} = \left\{ \beta_{i} : 1 \leq i \leq \ell \right\} $ and $\beta_{i} = \alpha_{i}^{\phi}$ with $\phi: B_{1} \to B_{2}$ the bijection discussed above.

    \begin{claim}
      If $\ell \geq 3$, the permutation group $G_{D_{1}}^{D_{1}}$ induced by the block stabiliser $G_{D_{1}}$  contains the subgroup $\left\langle \tau \right\rangle \times \sym(\ell) $ where $\tau$ is the superflip in $D_{1}$ swapping $\alpha_{i}$ and $\beta_{i}$ for $1 \leq i \leq \ell$.
    \end{claim}
    \begin{proof}\renewcommand{\qedsymbol}{}
      Since $\ell \geq 3$, we may assume that \cref{claim:Blarge} holds.
      First observe that since $G$ is transitive, there must be an element $g \in G$ that swaps $B_{1}$ and $B_{2}$.
      If $g$  preserves the sets $\left\{ \alpha_{i}, \beta_{i} \right\}$ for  $i \in \left\{ 1, \dots, \ell  \right\} $, then up to an element in $\left\{ 1 \right\} \times \sym(\ell)$, $g$ is the superflip. 
      Therefore, we may assume that $g$ does not preseve the set $\left\{ \alpha_{1}, \beta_{1} \right\}$.
      Moreover, without loss of generality we may assume that $\alpha_{1}^{g} = \beta_{1}$ and $\beta_{1}^{g} = \alpha_{2}$.
      Let $i$ and $j$ be such that $\alpha_{2}^{g} = \beta_{i} $ and $\beta_{2}^{g} = \alpha_{j}$.
      Observe that $i \neq 1$ and $j \neq 2$.
      This implies that $x^{g} = (\beta_{1} \ \beta_{i})(\alpha_{2}\ \alpha_{j}) \in G$, since $(\alpha_{1} \ \alpha_{i})( \beta_{1}\ \beta_{i}) \in G$, $(\alpha_{1}\ \alpha_{i})(\alpha_{2}\ \alpha_{j}) \in G$.
      Unless $i = 2$ and $j=1$, we have a contradiction.
      Similarly, if $r$ and $s$ are such that $\alpha_{3}^{g}=\beta_{r}$ and $\beta_{3}^{g} = \alpha_{s}$,
      then $( \alpha_{1}\ \alpha_{3})(\beta_{1}\ \beta_{3})^{g} = (\beta_{1} \ \beta_{r})(\alpha_{2}\ \alpha_{s}) \in G$ which implies that $(\alpha_{1}\ \alpha_{r}) (\alpha_{2}\ \alpha_{s}) \in G$ (observe that $r,s \geq 3$).
      It follows that there is either a $3$-cycle or a $2^{2}$-element (depending on whether or not $s=r$) whose support is contained in $B_{1}$, which contradicts our assumptions.
    \end{proof}
    
    \begin{claim}
      Assume that $\ell \geq 3$. For each $D \in \cD$, let $\sflip_{D} \in G $ be a permutation that induces the superflip in $G_{D}^{D}$.
      The family \[\cP := \bigcup_{D \in \cD} \left\{ \left\{ \alpha, \alpha ^{\sflip_{D}} \right\} : \alpha \in D \right\} \] is a block system for $G$.
    \end{claim}
    
    \begin{proof}\renewcommand{\qedsymbol}{}
      Assume otherwise and let $g \in G$ be an element that does not preserves $\cP$.
      We may assume, without loss of generality, that $g$ maps the pair $\left\{\alpha_{1} , \beta_{1} \right\} $ to a pair not in $\cP$.
      Label the elements of $D'= D^{g}_{1}$ such that $B_{1}^{g} = \left\{ \gamma_{i} : 1 \leq i \leq \ell \right\} $, $B_{2}^{g}= \left\{ \delta_{i}: 1 \leq i \leq \ell \right\} $ and $\delta_{i}= \gamma_{i}^{\sflip_{D'}}$.
      Moreover, we may assume that $\alpha_{1}^{g} = \gamma_{1}$ and $\beta_{1}^{g} = \delta_{2}$.
      Let $i \in \left\{ 1, \dots, \ell \right\} $  such that $\alpha^{g}_{i} = \gamma_{2}$;
      notice that $i \geq 2$.
      Observe that
      $(\alpha_{1}\ \alpha_{i})(\beta_{1} \ \beta_{i}) \in G$,
      which implies that
      $(\alpha_{1}\ \alpha_{i})(\beta_{1} \ \beta_{i})^{g}=(\gamma_{1} \ \gamma_{2})(\delta_{2}\ \delta_{j})  \in G$, where $1 \leq j \leq \ell$ is such that $\delta_{j} = \beta_{i}^{g}$.

      Hence $(\gamma_{1} \ \gamma_{2})(\delta_{1}\ \delta_{2})(\gamma_{1} \ \gamma_{2})(\delta_{2}\ \delta_{j}) = (\delta_{1}\ \delta_{2})(\delta_{2}\ \delta_{j}) \in G$.
      Unless $j = 1$ we have either a $3$ cycle or a $2^{2}$-element contained in $B_{2}^{g}$, which is a contradiction.
      Therefore $\beta_{i}^{g} = \delta_{1}$.
      In other words, we can conclude that
      \[\begin{aligned}
        \alpha_{1}^{g} &= \gamma_{1}, &
        \beta_{1}^{g} &= \delta_{2}, &
        \alpha_{i}^{g} &= \gamma_{2}, &
        \beta_{i}^{g} &= \delta_{1}.
      \end{aligned}\]
      
      Since $\ell\geq 3$, the set $\left\{ 1, \dots, \ell \right\} \sm \left\{ 1,i \right\}$ is not empty. 
      Let $t \in \left\{ 1, \dots, \ell \right\} \sm \left\{ 1,i \right\}$.
      Consider $r,s \in \left\{ 1, \dots, \ell \right\} $ such that $\gamma_{r} = \alpha_{t}^{g}$ and $\delta_{s} = \beta_{t}^{g}$ and observe that $\left\{ r,s\right\} \cap \left\{ 1,i \right\} = \emptyset$.
      Finally observe that
      \[
        (\alpha_{1}\ \alpha_{t})(\beta_{1}\ \beta_{t})^{g} = (\gamma_{1}\ \gamma_{r})(\delta_{2}\ \delta_{s}) \in G
        \]
        which implies that $(\gamma_{1}\ \gamma_{r})(\delta_{1}\ \delta_{r})(\gamma_{1}\ \gamma_{r})(\delta_{2}\ \delta_{s}) = (\delta_{1}\ \delta_{r})(\delta_{2}\ \delta_{s}) \in G$, but this element is a $3$-cycle (if $s=r$) or a $2^{2}$-element (if $s \neq r$) contained in $B_{2}^{g}$, which as before, leads to a contradiction.
      \end{proof}
      
      Observe that the block system $\cP$ defined above consists of minimal  blocks of size $2$.
      Notice also that the $2^{2}$-element $(\alpha_{1} \ \alpha_{2})(\beta_{1} \ \beta_{2})$ swaps two blocks in $\cP$.
      It follows that if $\ell \geq 3$, \cref{item:4easy} holds, that is, $G$ is one of the groups described in \cref{prop:Easy}.

      On the other hand, if the blocks of $\cB$ are of size $2$, so that \cref{claim:Blarge} does not hold, then each block of $\cD$ contains in general $m \geq 2$ blocks of $\cB$.
      In this case, the group $G_{D}^{D}$ is permutation isomorphic to a subgroup of $\cubeWr$ and the group $G_{(\Omega \sm D)}$ contains a group $X$ isomorphic to $\evenAb[m]$  generated by the conjugates of $x$ under $G_{D}$. 
      This is of course equivalent to the fact that \cref{item:4hard} holds.
    \end{proof}

\cref{prop:Hard} provides us with the last tool to complete the proof of \cref{thm:groups_minDeg4}.

\begin{proof}[Proof of \cref{thm:groups_minDeg4}]
  Let $G$ be a transitive permutation group containing a $2^{2}$-element $x=(\alpha_{1} \alpha_{2})(\beta_{1}\ \beta_{2} )$. Let $\cB$ be a system of minimal blocks and let $B_{1}$ be block of $\cB$ that contains $\alpha_{1}$.
  \cref{rem:deg4_supp} claims that we have three cases:

  Asuume that $\supp(x) \subset B_{1}$. 
  In this case, as a direct consequence of \cref{lem:theOne}, \cref{lem:DTSmallDeg} and \cref{lem:babai} and an easy computation in GAP \cite{TheGAPGroup_2022_GapGroupsAlgorithmsGapGroupsAlgorithms}, we conclude that $G$ is as in \cref{item:groups_minDeg4_prim} of \cref{thm:groups_minDeg4}.

  Assume now that $\supp(x) \cap B_{1} = \left\{ \alpha_{1}, \beta_{1} \right\}$. In this situation $G$ satisfies \cref{prop:Easy}, which directly implies that \cref{item:groups_minDeg4_cross} of \cref{thm:groups_minDeg4} holds.

  Finally, if $\supp(x) \cap B_{1} = \left\{ \alpha_{1}, \alpha_{1} \right\}$ the group $G$ must satisfy \cref{prop:Hard}. 
  If any of \cref{item:lessThan4} or \ref{item:4prim} of \cref{prop:Hard} holds, then $G$ satisfies \cref{item:groups_minDeg4_prim} of \cref{thm:groups_minDeg4}; on the other hand if \cref{item:4easy} or \cref{item:4hard} holds, then $G$ satisfies \cref{item:groups_minDeg4_cross} of \cref{thm:groups_minDeg4}.
\end{proof}

 \section{Vertex transitive graphs with small motion} \label{sec:graphs}
Recall that the motion $\mu(\Gamma)$ of a vertex-transitive graph $\Gamma$ is the minimal degree of its automorphism group.
In this section we use the results obtained in \cref{sec:groups} to classify the vertex-transitive graphs whose motion is $4$ or a prime number $p$.
Our aim is to prove \cref{thm:graphs}, which we repeat below:

\begin{thm*}
  Let $\Gamma$ be a vertex-transitive graph on $n$ vertices. 
Then:
  \begin{enumerate}
    \item \label{item:graphs_mindeg2} 
    If $\aut(\Gamma)$ contains a $2$-cycle (or equivalently if $\mu(\Gamma) = 2$) then
$\Gamma \cong \lex{K_{m}}{\Theta}$ or $\Gamma \cong \lex{(mK_{1})}{\Theta}$ with
    $\Theta$ a vertex-transitive graph on $k$ vertices and $m \geq 2$ such that $n =mk$.  
    Conversely, for every vertex-transitive graph $\Theta$ and $m \geq 2$ the graphs $\lex{K_{m}}{\Theta}$ and $\lex{(mK_{1})}{\Theta}$ have motion $2$.
    \item \label{item:graphs_mindegp} If $\aut(\Gamma)$ contains a $p$-cycle for some prime number $p \geq 3$, then either:
    \begin{enumerate}
      \item $\mu(\Gamma) = 2$ and thus $\Gamma$ is one of the graphs described in (\ref{item:graphs_mindeg2}); or
      \item $\Gamma$ is isomorphic to $\lex{\Sigma_{p}}{\Theta}$ with $\Sigma_{p}$ a circulant graph with $p$ vertices and $\Theta$ is a vertex-transitive graph; in this case $\mu(\Gamma) = 2$ if $\Sigma_{p}$ is isomorphic to $ K_{p} $ or $pK_{1}$, and $\mu(\Gamma) = p-1$ otherwise. 
    \end{enumerate}
\item \label{item:graphs_mindeg4} If $\mu(\Gamma) = 4$ then one of the following holds:
      \begin{enumerate}
        \item 
        \label{item:graphs_mindeg4_prim}
        $\Gamma \cong \lex{C_{5}}{\Theta}$ with $\Theta$ a vertex-transitive graph; or \item \label{item:graphs_mindeg4_cross1}
         $\Gamma \cong \lex{(K_{m} \square K_{2})}{\Theta}$ or $\Gamma \cong \lex{\bar{K_{m} \square K_{2}}}{\Theta}$ with $m\geq 3  $ and $\Theta$ a vertex-transitive graph; or
\item \label{item:graphs_mindeg4_cross2} $\Gamma \cong  \Inf$, defined in \cref{def:inf} below, for some 
        $m \geq 2$, $\lambda,\kappa \in \bZ_{2}$, a graph $\Sigma$ and a partition $\cP$ of $V\Sigma$ with blocks of size $2$ satisfying some natural conditions precised in \cref{lem:inf}.
        \end{enumerate}
      Conversely, for every vertex-transitive graph $\Theta$, every vertex-transitive graph $\Sigma$ admitting a $\aut(\Sigma)$-invariant partition $\cP$ of its vertices with blocks of size $2$ satisfying the conditions in \cref{lem:inf}, and  $m \geq 2$, $\lambda, \kappa \in \bZ_{2}$, the graphs $\lex{C_5}{\Theta}$, $\lex{(K_{m} \square K_2)}{\Theta}$, $\lex{\bar{K_{m} \square K_{2}}}{\Theta}$ and $\Inf$  have motion $4$.
  \end{enumerate}
\end{thm*}

The proof will be carried over in a series of results. 
First we shall classify the vertex transitive graphs with motion $2$ (see \cref{thm:graphsMinDeg2}), showing that \cref{item:graphs_mindeg2} of \cref{thm:graphs} holds.
Then we shall use \cref{thm:groups_minDegp} to show that \cref{item:graphs_mindegp} of \cref{thm:graphs} holds; in particular we shall prove that there are no vertex transitive graphs with motion a prime number $p \geq 3$ (see \cref{coro:graphsMinDeg3,coro:graphsMinDegP}).
Finally, we shall split the classification of the vertex transitive graphs of motion $4$ using the results in \cref{sec:grpsMinDeg4} (see \cref{lem:graphs4}).
In \cref{thm:graphs_prim} we show that if \cref{item:graphs4prim} of \cref{lem:graphs4} holds, then a vertex transitive graph with motion $4$ must be isomorphic to $\lex{C_{5}}{\Theta}$ with $\Theta$ a vertex transitive graph (\cref{item:graphs_mindeg4_prim} of \cref{thm:graphs}).
Then we turn our analysis to the case where $\Gamma$ satisfies \cref{item:graphs4easy} of \cref{lem:graphs4}, and show that either $\mu(\Gamma) < 4$ or $\Gamma$ satisfies either \cref{item:graphs_mindeg4_cross1} (\cref{thm:graphs4_lex}) or \cref{item:graphs_mindeg4_cross2} (\cref{prop:graphs_inf}) of \cref{thm:graphs}.
Ultimately, we prove that if $\Gamma$ satisfies \cref{item:graphs4hard} of \cref{lem:graphs4} then $\mu(\Gamma)=2$ or $\Gamma$ is as described in \cref{item:graphs_mindeg4_cross2} of \cref{thm:graphs} (\cref{prop:graphs_hard}), thus completing the proof of \cref{thm:graphs}.

To begin with, we establish some classical graph theoretical notation.
If $\Gamma$ is a graph with vertex set $V\Gamma$ and $B \subset V\Gamma$, the \emph{subgraph induced by $B$} is the graph $\Gamma[B]$ whose vertex set is $B$ and such that $\beta_{1} \sim \beta_{2}$ in $\Gamma[B]$ if and only if $\beta_{1} \sim \beta_{2}$ in $\Gamma$.
Informally speaking, $\Gamma[B]$ is the graph resulting after deleting all the vertices not in $B$ and all the edges incident to any of those vertices.

If $\alpha$ is a vertex of $\Gamma$, then $\Gamma(\alpha)$ denotes the \emph{neighbourhood} of $\alpha$, that is, the set \[\Gamma(\alpha) = \left\{ \beta \in V\Gamma : \alpha \sim \beta \right\}. \]

Assume that $\cB$ is a partition of $V\Gamma$. Then the graph $\Gamma/\cB$ is the graph whose vertices are the sets in $\cB$ and such that $B_{1} \sim B_{2}$ in $\Gamma / \cB$ if there exist vertices $\alpha_{1} \in B_{1}$ and $\alpha_{2} \in B_{2}$ such that $\alpha_{1} \sim \alpha_{2}$.
The graph $\Gamma/\cB$ is called the \emph{quotient graph}.
Of our particular interest is when $\cB$ is a block system for a group  $G \leq \aut(\Gamma)$;
 in the particular case when $\cB$ is the set of orbits of $G$ the graph $\Gamma/ \cB$ is also denoted by $\Gamma/G$.

If $\Gamma$ and $\Delta$ are graphs with vertex set $V\Gamma$ and $V\Delta$ respectively, the \emph{lexicographic product} of $\Gamma$ and $\Delta$ is the graph $\lex{\Delta}{\Gamma}$ (also denoted by $\Gamma[\Delta]$ by some authors) defined by
\[\begin{aligned}
	V(\lex{\Delta}{\Gamma}) &= V\Delta \times  V\Gamma \\
	(\delta_{1},\gamma_{1}) \sim (\delta_{2}, \gamma_{2})
&\iff
	\begin{cases}
	    \delta_{1} \sim \delta_{2} \text{ and } \gamma_{1} = \gamma_{2} & \text{or }\\
	  \gamma_{1} \sim \gamma_{2}
	\end{cases}
\end{aligned}\]

The following results are elementary but they shall be very useful. Their proofs are straightforward and are omitted.

\begin{lem}\label{lem:graphsTrans}
  Let $\Gamma$ be graph with vertex set $V\Gamma$ and
  let $\omega_{1}, \omega_{2} \in V\Gamma$.
  The transposition $x = (\omega_{1} \ \omega_{2})$ is an automorphism of $\Gamma$ if and only if $\Gamma(\omega_{1}) \sm \left\{ \omega_{2} \right\}  = \Gamma(\omega_{2}) \sm \left\{ \omega_{1} \right\} $.
\end{lem}

\begin{lem}\label{lem:lexicographic}
  Let $\Gamma$ be a graph with vertex set $V\Gamma$  and $G \leq \aut(\Gamma)$.
  Assume that $G$ acts imprimitively on $V\Gamma$ and let $\cB$ be a block system for $G$.
  Assume also that for some $B \in \cB$ the group $G_{(V\Gamma \sm B)}$ is transitive on $B$.
  Then $\Gamma$ is isomorphic to the  graph $\lex{\Gamma[B]}{(\Gamma/\cB)}$.
\end{lem}

We are now ready to classify vertex-transitive graphs whose motion is a prime number.

\begin{thm}\label{thm:graphsMinDeg2}
  Let $\Gamma$ be a vertex-transitive graph whose motion is $2$.
  Then $\Gamma \cong \lex{K_{m}}{\Theta}$ or $\Gamma \cong \lex{(mK_{1})}{\Theta}$ for some $m \geq 2$ and a vertex-transitive graph $\Theta$.
	Moreover, all such graphs have motion $2$.
\end{thm}
\begin{proof}
  This is a direct consequence of \cref{thm:groups_minDegp} and \cref{lem:lexicographic}.
  If $G=\aut(\Gamma)$ is of minimal degree $2$ then we may assume that $V\Gamma = [m] \times [k]$ and think of  $G$ as subgroup of $  \sym(m) \wr \sym(k)$ preserving the partition $\cB = \left\{ B_{i}  : 1 \leq i \leq k \right\}$, where $B_{i}= [m] \times \left\{ i \right\}$ and such that the base group $\sym(m)^{k}$  is contained in $ G$.
  Since $\sym(m)$ is transitive on $[m]$, \cref{lem:lexicographic} implies that $\Gamma \cong \lex{\Delta}{\Theta}$ with $\Theta = \Gamma/\cB$ and $\Delta =\Gamma[B_{1}]$.
  However, $\Delta$ must be isomorphic to $K_{m}$ or to $mK_{1}$, because $\sym(m)$ is doubly transitive on $[m]$.
  The graph $\Gamma/\cB$ is vertex-transitive because $G$ is transitive on $\cB$.

  Conversely, let $\Theta$ be a vertex-transitive graph and $m \geq 2$.
  It is known that $(\aut(\Delta) \wr \aut(\Theta)) \pleq \aut(\lex{\Delta}{\Theta})$. 
  It follows that if $\Delta$ and $\Theta$ are vertex-transitive graphs, then so it is $\lex{\Delta}{\Theta}$. 
  Moreover, if $\Delta \cong K_{m}$ or $\Delta \cong m K_{1} $, then $\aut(\lex{\Delta}{\Theta})$ contains a transposition, which implies that $\mu(\Gamma) = 2$.
\end{proof}

\begin{coro}\label{coro:graphsMinDeg3}
  There are no vertex-transitive graphs whose motion is $3$. 
  More precisely, if the automorphism group of a graph contains a $3$-cycle $x$, then its automorphism group is of minimal degree $2$.
\end{coro}
\begin{proof}
  The result follows from \cref{thm:groups_minDegp} and a similar argument to that in \cref{thm:graphsMinDeg2}. 
  Recall that $\alt(m)$ is transitive on unordered pairs if $m \geq 3$.
  We can prove that if $\aut(\Gamma)$ contains a $3$-cycle, then $\Gamma  \cong \lex{K_{m}}{\Theta}$ or $ \Gamma \cong \lex{(mK_{1})}{\Theta}$ for some vertex-transitive $\Theta$ and $m \geq 3$.
  However, according to \cref{thm:graphsMinDeg2}, these graphs have motion $2$.
\end{proof}

\cref{coro:graphsMinDegP} is a direct consequence of \cref{thm:groups_minDegp}, which in turn follows from the results in \cite{Jones_2014_PrimitivePermutationGroups} (see \cref{thm:primitiveP}). 
However, as mentioned in \cref*{sec:intro}, the results in \cite*{Jones_2014_PrimitivePermutationGroups} rely on the classification of finite simple groups (CFSG). 
The following results allow us to obtain a proof for \cref{coro:graphsMinDegP} that does not rely on the CFSG.

The following result is a theorem often credited to Marggraff \cite{Marggraff1889_UeberPrimitiveGruppen_PhDThesis,LevinTaylo1976_TheoremMarggraffPrimitive} but originally proved by Jordan \cite{Jordan_1871_TheoremesSurLes,Neumann1985_SomePrimitivePermutation}. 
See \cite[Section 4]{Jones_2014_PrimitivePermutationGroups} for a historical note on the topic.

\begin{thm}\label{thm:JordanTr}
  Let $G$ be a primitive permutation group on a finite set $\Omega$. 
  Suppose that $G$ contains a non-trivial cyclic subgroup which fixes $k$ points of $\Omega$ and which is transitive on the remaining points. 
  In this situation $G$ is $(k+1)$-transitive.
\end{thm}

Naturally, following an analogous analysis to that in \cref{sec:groups}, we obtain the following result.

\begin{coro}\label{coro:JordanNP}
  Let $G$ be permutation group of a set $\Omega$ with $|\Omega| = n$. 
  Assume that for some prime number $p$, $G$ contains a $p$ cycle. 
  Then there exist $m,k \in \bN$ with $p \leq m$ and $n = mk$,  and transitive permutation groups $X \leq Y$ of degree $m$ such that $X$ is primitive and contains a $p$-cycle,  and 
  \[X^{|k|} \leq G \leq Y \wr \sym(k).\] 
  Moreover, if $p < m$, then $X$ is $2$-transitive.
\end{coro}

\begin{coro}\label{coro:graphsMinDegP} 
  There are no vertex-transitive graphs with motion $p$ for prime number $p\geq 5$. 
  More precisely, if $\Gamma$ is a vertex-transitive graph whose automorphism group contains a $p$-cycle, then $\Gamma \cong \lex{\Delta}{\Theta}$ where $\Delta$ is either $K_m$ or $m K_1$ for some $m \geq 2$ or a circulant graph $\Sigma_p$ with $p$ vertices, and $\Theta$ is a vertex-transitive graph.
\end{coro}
\begin{proof}
  This is a direct consequence of \cref{coro:JordanNP}.
  Assume that $\Gamma$ is a graph whose automorphism group $G$ contains a $p$-cycle.
  Then \[X^{|k|} \leq G \leq Y \wr \sym(k)\] for some transitive permutation group $X \leq Y $ of degree $m$ such that $|V\Gamma|= mk$, $p \leq m$ and $X$ is primitive and contains a $p$-cycle.
  We may assume that $V\Gamma = [m] \times [k]$ (with $X$ and $Y$ acting on $[m]$).
  Since $X$ is transitive, \cref{lem:lexicographic} then implies that $\Gamma \cong \lex{\Delta}{(\Gamma/\cB)}$ with $\Delta \cong \Gamma([m] \times \left\{ 1 \right\})$.
  If $m = p $, the graph $\Delta$ admits a cyclic subgroup acting regularly on the vertices, hence it is isomorphic to a circulant graph $\Sigma_{p}$. It is known that if $\Sigma_{p}$ is a circulant graphs with vertex set $C_{p}$, the mapping $z \mapsto z^{-1}$ is an automorphism of $\Sigma_{p}$ fixing at least one vertex, and since $p\geq 5$, then this automorphism is non trivial. 
  It follows that $\aut(\Delta)$ contains an automorphism $x$ with $\supp(x) < p$, which implies that $\mu(\Gamma) \leq \mu(\Delta) < p$.

  On the other hand if $p < m$,  the group $X$ is $2$-transitive, which implies that $\Delta$ is either complete or edgeless; which in turn implies that $\Gamma$ is as in \cref{thm:graphsMinDeg2}, that is, $\mu(\Gamma) = 2$.     
  \end{proof}

We shall divide the classification of vertex-transitive graphs of motion $4$ in a similar fashion as we discussed transitive groups of minimal degree $4$ in \cref{sec:grpsMinDeg4}.
More precisely, let us consider the following (obvious) lemma.
\begin{lemma} \label{lem:graphs4}
  Let $\Gamma$ be a vertex transitive graph and $G$ a vertex-transitive group of automorphisms of $\Gamma$ satisfying that $\mu(G) = 4$. Then 
  \begin{enumerate}
    \item \label{item:graphs4prim} $G$ is as in \cref{item:groups_minDeg4_prim} of \cref{thm:groups_minDeg4} and we may assume that $V\Gamma = [m] \times [k]$ and \[X^{k} \leq G \leq Y \wr \sym(k) \] with $m,X,Y$ as in \cref{tab:22-prim}; or
\item \label{item:graphs4easy} $G$ is as in \cref{prop:Easy}, and we may assume that $V\Gamma = \left\{ 1, -1, 2,-2, \dots,  m, -m \right\} \times [k]$, and \[X^{k} \leq G \leq Y \wr \sym(k)\] with $(X,Y)$ as in \cref{tab:22-easy}; or
    \item \label{item:graphs4hard} $G$ is as in \cref{item:4hard} of \cref{prop:Hard}, and we may assume that $V\Gamma = \left\{ 1, -1, 2,-2, \dots,  m, -m \right\} \times [k]$ and \[(\evenAb)^{k} \leq G \leq Y \wr \sym(k)\] for some $Y \leq \cubeWr$.
  \end{enumerate}
\end{lemma}

Let us assume that $G$ is as in \cref{item:graphs4prim} of \cref{lem:graphs4} above.
Unless $m=5$ and $(X,Y)= (\dih(5), \agl_{1}(5))$,  the group $X$ is doubly transitive.
Just as in \cref{thm:graphsMinDeg2}, this implies that $\Gamma \cong \lex{K_m}{(\Gamma/\cB)}$ or $\Gamma \cong  \lex{mK_{1}}{(\Gamma/\cB)}$, with $\cB = \left\{ [m] \times \left\{ i \right\} : 1 \leq i \leq k \right\}$. 
As shown before, these graphs have motion $2$.

On the other hand, if $(m,X,Y) = (5, \dih(5), \agl_{1}(5))$, then  \cref{lem:lexicographic} implies that $\Gamma\cong \lex{\Delta}{(\Gamma/\cB)}$ where $\Delta = \Gamma\left[ [m] \times \left\{ 1 \right\} \right]$ is a graph preserved by $\dih(5)$.
These graphs are in correspondence with a subset of the orbits of $\dih(5)$ on the edges of $K_{5}$.
There are just two such orbits each of them leading to isomorphic graphs $C_{5}$.
It follows that $\Delta$ must be complete, edgeless or $C_{5}$ and as seen before, in the first two cases $\Gamma$ has motion $2$.
\cref{lem:graphsTrans} proves that $\lex{C_{5}}{\Theta}$ does not have minimal degree $2$ for any vertex-transitive graph $\Theta$ since there are no two vertices $\omega_{1}, \omega_{2}$ satisfying $\Gamma(\omega_{1}) \sm \left\{ \omega_{2} \right\}  = \Gamma(\omega_{2}) \sm \left\{ \omega_{1} \right\} $. 
The discussion above is summarized in the following result.

\begin{thm}\label{thm:graphs_prim}
  If $\Gamma$ is a vertex transitive graph with vertex set $V\Gamma= [m] \times [k]$ and $G \leq \aut(\Gamma)$ satisfies that \[ X^{k} \leq G \leq Y \wr \sym(k) \] for some $X,Y$ as in \cref{tab:22-prim} then 
  \begin{enumerate}
    \item $\mu(\Gamma)=2$ and $\Gamma \cong \lex{K_{m}}{\Theta}$ or $\Gamma \cong \lex{(mK_{1})}{\Theta}$ for a vertex transitive graph $\Theta$; or
    \item $\mu(\Gamma) = 4$ and $\Gamma\cong \lex{C_{5}}{\Theta}$ for a vertex transitive graph $\Theta$.
  \end{enumerate}
\end{thm}

Now  assume that \cref{item:graphs4easy} of \cref{lem:graphs4} holds,
meaning that we may assume that $V\Gamma = \left\{ 1, -1, 2,-2, \dots,  m, -m \right\} \times [k]$ for some $m,k$ such that $|V\Gamma| = 2mk$ and $G$ satisfies \[ X^{k} \leq G \leq Y \wr \sym(k) \] for a pair of groups $(X,Y)$ as in \cref{tab:22-cross}.

Let $D = \left\{ 1, -1, 2,-2, \dots,  m, -m \right\} \times \left\{ 1 \right\} $ and consider the group $Z = G_{(V\Gamma \sm  D)}$, that is, the pointwise stabiliser of the complement of the block  $D$.
This group satisfies $X \leq Z \leq Y$. 
If the group $Z$ is transitive, then \cref{lem:lexicographic} implies that $\Gamma \cong \lex{\Gamma[D]}{(\Gamma/\cD)}$ where  $\cD = \left\{ \left\{ 1, -1, 2,-2, \dots,  m, -m \right\} \times \left\{ i \right\} : 1 \leq i \leq k \right\}$.
For every pair $(X,Y)$ in \cref{tab:22-cross}, $X$ is an index-$2$ subgroup of $Y$ and the only case where $Z$ is not transitive on $D$ is when \[ \left\{ 1 \right\} \times \sym(m) = X = Z < Y = \left\langle \sflip \right\rangle \times \sym(m).  \] 
 We shall leave this case for later and for now we classify the graphs $\Delta$ with vertex set $\left\{ 1, -1, 2, -2, \dots, m,-m \right\}$ preserved by $Z$ when $Z = \left\langle \sflip \right\rangle  \times \sym(m)$ or $\evenAb \rtimes \sym(m) \leq Z$.
The following lemma describes such graphs when $Z = \left\langle \sflip \right\rangle  \times \sym(m)$.

\begin{lem} \label{lem:graphs_easySF}
  Let $m \geq 2$ and let $\Delta$ be a graph with vertex set $V\Delta= \left\{ 1, -1, 2, -2, \dots, m,-m \right\}$.
  Assume that $\Delta$ is preserved by the group $Z= \left\langle \sflip \right\rangle \times \sym(m)$.
  Then $\Delta$ is isomorphic to $K_{m} \square K_{2}$, $K_{m,m}$, $mK_{2}$, $K_{2m}$, or the complement of one of these.
  Moreover, if $\mu(\Delta) = 4$, then $m \geq 3$ and $\Delta$ is isomorphic to either $K_{m} \square K_{2}$ or $\bar{K_{m} \square K_{2}}$.
\end{lem}
\begin{proof}
  As before, any graph preserved by the group $Z$ is in correspondence with a family of orbits of $Z$ on unordered pairs of $V\Delta$.
  The orbits on pairs of $V\Delta$ of the group $Z$ are the following:
  \begin{align*}
    O_{1} &= \left\{ \left\{ i,j \right\} : 0 < i \neq j \leq m\right\} \cup \left\{ \left\{ i,j \right\} : -m \leq i \neq j < 0\right\},\\
    O_{2} &= \left\{ \left\{ i,-i \right\} : 1 \leq i \leq m \right\},\\
    O_{3} &= \left\{ \left\{ i,-j \right\} :1 \leq i \neq j \leq m \right\}.
  \end{align*} 
  The families $\left\{ O_{1}, O_{2} \right\}$, $\left\{ O_{2}, O_{3} \right\}$, $\left\{ O_{2} \right\}$ and $\left\{ O_{1}, O_{2}, O_{3} \right\}$ give rise to the graphs $K_{m} \square K_{2}$, $K_{m,m} $, $m K_{2}$ and $K_{2m}$, respectively; while their complementary families induce the corresponding complementary graphs. 
  Clearly the graphs $K_{m,m}$, $mK_{2}$, and $K_{2m}$ have motion $2$, while the graph $K_{m} \square K_{2}$ has a pair a of vertices satisfying \cref{lem:graphsTrans} if and only if $m=2$, which in turn implies that if $\mu(\Delta) = 4$, then $m\geq3$ and $\Delta$ is either $K_{m} \square K_{2}$ or its complement.
  Finally observe that for  $m \geq 3$,   $2 < \mu(K_{m} \square K_{2}) \leq 4 $;
  however, \cref{coro:graphsMinDeg3} implies that $\mu(K_{m} \square K_{2})$ cannot be $3$.
\end{proof}

The following result is the analogous to \cref{lem:graphs_easySF} but for the case when $\evenAb \times \sym(m) \leq Z$.
We omit the proof because it follows exactly the same ideas of that of \cref{lem:graphs_easySF}.

\begin{lem}\label{lem:graphs_easyTrans}
  Let $m \geq 3 $ and $\Delta$ be a graph with vertex set $V\Delta= \left\{ 1, -1, 2, -2, \dots, m,-m \right\}$.
  Assume that $\Delta$ is preserved by the group $\evenAb \times \sym(m)$.
  Then $\Delta$ is isomorphic to $K_{2m}$, $\bar{K_{2m}}$, $m K_{2} \cong \lex{K_{2}}{mK_{1}}$ or $\bar{m K_{2}} \cong \lex{2K_{1}}{K_{m}}$. 
  In particular, all such graphs $\Delta$ have motion $2$.
\end{lem}

Notice that if $m=2$, then $\evenAb \rtimes \sym(m) = \left\langle \sflip \right\rangle \times \sym(m) $.
In this case \cref{lem:graphs_easySF} applies and the resulting graphs are $K_{2} \square K_{2} \cong \bar{2 K_{2}}$, $K_{2,2} \cong \bar{2 K_{2}}$, $2 K_{2}$, $K_{4}$ or their complements. 
This shows that the conclusion of \cref{lem:graphs_easyTrans} holds even for $m=2$. 

We are now ready to continue with the classification of vertex-transitive graphs $\Gamma$ whose automorphism group contains a group $G$ satisfying \cref{item:graphs4easy} of \cref{lem:graphs4}.
As mentioned before, if $D = \left\{ 1, -1, 2,-2, \dots,  m, -m \right\} \times \left\{ 1 \right\} $ and  $Z = G_{(V\Gamma \sm  D)}$ is transitive on $D$, then $\Gamma \cong \lex{\Delta}{(\Gamma/\cD)}$ with $\Delta$ a graph preserved by $Z$.  
The possibilities for $\Delta$ are discussed in \cref{lem:graphs_easySF} and \cref{lem:graphs_easyTrans}; but observe that all such graphs, except $K_{m} \square K_{2}$ and $\bar{K_{m} \square K_{2}}$ (for $m \geq 3$) admit a transposition as an automorphism.
This transposition then induces  a transposition in $\Gamma \cong \lex{\Delta}{(\Gamma/\cD)}$.
This shows that unless $\Delta \cong K_{m} \square K_{2}$ or $\Delta \cong \bar{K_{m} \square K_{2}}$ with $m\geq 3$, the motion of $\Gamma$ is $\mu(\lex{\Delta}{(\Gamma/\cD)}) = 2$.
In summary,

\begin{thm}\label{thm:graphs4_lex}
Let $\Gamma$ be a vertex transitive graph with vertex set $V\Gamma = \left\{ 1, -1, 2,-2, \dots,  m, -m \right\} \times [k]$. Assume that $\Gamma$ admits a group $G \leq \aut(\Gamma)$ such that 
\[X^{k} \leq G \leq Y \wr \sym(k)\] with $(X,Y)$ as in \cref{tab:22-easy}, in such a way that $G_{(V\Gamma\sm \left\{ 1, -1, 2,-2, \dots,  m, -m \right\} \times \left\{ 1 \right\} })$ is transitive on $\left\{ 1, -1, 2,-2, \dots,  m, -m \right\} \times \left\{ 1 \right\} $. 
Then 
\begin{itemize}
  \item $m \geq 3$ and $\Gamma\cong \lex{\Delta}{\Theta}$ with $\Theta$ a vertex transitive graph and $\Delta \in \left\{ K_{m} \square K_{2}, \bar{K_{m} \square K_{2}} \right\}$. In this case, $\mu(\Gamma) = 4$; or
  \item $\Gamma\cong \lex{\Theta}{\Delta}$ with $\Delta \in \left\{ K_{2m}, \bar{K_{2m}}, mK_{2}, \bar{mK_{2}} \right\}$. In this case, $\mu(\Gamma) = \mu(\Delta)=2$.
\end{itemize}
\end{thm}
 
The previous analysis covers all the possible vertex transitive graphs with motion $4$ that are isomorphic to lexicographic products. 
To complete the classification of such graphs we shall use the construction $\Inf$, which was introduced in \cref{def:inf} and that we revisit below.
Let $\Sigma$ be a graph with vertex set $V\Sigma$ and let $\cP$ be a partition of $V\Sigma$ with sets of size $2$.
Let $m \geq 2$ be an integer and let $\lambda,\kappa \in \bZ_{2} $.
The graph $\Inf$ is the graph whose vertex set is $V\Sigma \times \left\{ 1, \dots, m \right\} $ and the edges are given by
\[
(\alpha, i) \sim (\beta, j) \iff
\begin{cases}
  \left\{ \alpha, \beta \right\} \not\in \cP
  \text{ and }
  \alpha \sim \beta \text{ in } \Sigma
  \\
  \left\{ \alpha,\beta \right\} \in \cP,
  \lambda = 1,
  \text{ and }
  i=j, \\
  \left\{ \alpha,\beta \right\} \in \cP,
  \lambda = 0,
  \text{ and }
  i\neq j,\\
  \alpha = \beta \text{ and } \kappa = 1.
\end{cases}
\]

Now we shall establish some key symmetry properties of the construction $\Inf$.

\begin{lem}\label{lem:inf}
  Let $\Sigma$ be a graph with vertex-set $V\Sigma$ and let $\cP$ be a partition of $V\Sigma$ with sets of size $2$. 
  Let $m \in \bN$ and $\lambda, \kappa \in \bZ_{2}$.
  Let $\aut_{\cP}(\Sigma)$ denote the group of automorphisms of $\Sigma$ that preserve $\cP$, that is \[\aut_{\cP}(\Sigma) = \left\{ x \in \aut(\Sigma) : \left\{ \alpha,\beta \right\} \in \cP \iff \left\{ \alpha^{x},\beta^{x} \right\} \in \cP \right\}.\] Then the following hold:
  \begin{enumerate}
    \item\label{item:infVT} The graph $\Inf$ is vertex-transitive if and only if the group $\aut_{\cP}(\Sigma)$ is transitive on $V\Sigma$.
  \item \label{item:mInf4} If $\Inf$ is vertex-transitive, then $\mu(\Inf) = 4$ unless $\kappa = \lambda$ and there is a block $\left\{ \alpha, \beta\right\} \in \cP$ such that the permutation $(\alpha\ \beta)$ is an automorphism of $\Sigma$.
  \end{enumerate}
\end{lem}
\begin{proof}
One can verify that every automorphism of $\Inf$ projects to an automorphism of $\Sigma$ that preserves $\cP$.  
Conversely, if $x \in \aut_{\cP}(\Sigma)$, then the mapping $\bar{x}: V\Sigma \times \left\{ 1, \dots, m \right\} \to V\Sigma \times \left\{ 1, \dots, m \right\}  $ given by $(\alpha,i)^{\bar{x}} = (\alpha^{x}, i)$ is an automorphism of $\Inf$.
For every pair $\left\{ \alpha, \beta \right\} \in \cP$, the group $\sym(m)$ acts on the set $\{ \alpha, \beta\} \times [m] $ by permuting the second coordinates.
This action induces in turn a subgroup of $\aut(\Inf)$ isomorphic to $\sym(m)^{|\cP|}$ that preserves the family
\[\bar{\cP} = \left\{ \left\{ \alpha, \beta \right\} \times [m] :\left\{ \alpha, \beta \right\} \in \cP  \right\} \] setwise.
The discussion above implies that $\Inf$ is vertex transitive if and only if the group $\aut_{\cP}(\Sigma)$ acts transitively on $V\Sigma$.

Observe that if $\left\{ \alpha,\beta \right\} \in \cP $, then the permutation $\left( (\alpha, 1) \ (\alpha, 2) \right) ( (\beta,1)\ (\beta,{2}) )$ is an element of $ \aut(\Inf)$.
This implies that $\Inf$ has motion at most $4$.
An easy consequence of \cref{lem:graphsTrans} is that $\Inf$ has motion $2$ if and only if $\lambda \neq  \kappa$ and $\Sigma$ admits a transposition preserving $\cP$.
In this case the family
 \[
   \big\{ \left\{ (\alpha, i), (\beta, i) \right\}  : \left\{ \alpha,\beta \right\} \in \cP \text{ and } i \in \left\{ 1, \dots, m \right\} \big\}
 \]
  is a block system for $\aut(\Inf)$ with blocks of size $2$ such that the transposition swapping the two elements of a given block is an automorphism of $\Inf$.
\end{proof}

  Observe that if $\Sigma$ admits a transposition $(\alpha \ \beta)$ for some $\left\{ \alpha, \beta \right\} \in \cP$, but $\lambda = \kappa$, then $\Inf \cong (\Sigma/\cP)[\Delta]$ with $\Delta = K_{m} \square K_{2}$ (if $\lambda = \kappa=1$) or $\Delta = \bar{K_{m} \square K_{2}}$ (if $\mu = \kappa =0$) which we already proved that have motion $4$.

We shall use the construction $\Inf$ described above to complete the classification of vertex-transitive of motion $4$.

It remains to explore the vertex transitive graphs $\Gamma$ with vertex set $ V\Gamma = \left\{ 1, -1, 2, -2, \dots, m,-m \right\} \times [k]$ such that if we let $X$ and $Y$ denote the groups $\left\{ 1 \right\} \times \sym(m) $ and $ \left\langle \sflip \right\rangle \times \sym(m) $, respectively then $\aut(\Gamma)$ contains a group $G$ satisfying:
\begin{equation}\label{eq:group_inf}
  X^{k} \leq G \leq Y \wr \sym(k) 
\end{equation}
Moreover, in this case we might actually assume that $\aut(\Gamma) = G$; otherwise $\mu(\Gamma) < 4$ or $\mu(\Gamma)=4$ and $\Gamma$ is a lexicographic product as described before. 
If we let $D = \left\{ 1, -1, 2, -2, \dots, m,-m \right\} \times \left\{ 1 \right\}$, we might also assume that the pointwise stabiliser $G_{(V\Gamma \sm D)}$ of the complement of $D$ is exactly $X$. 

We shall prove that in this situation $\Gamma \cong \Inf$ for a suitable choice of  $\Sigma$ and $\cP$. 
First observe that the vertex set of $\Gamma$ already possesses a structure that resembles that of a graph $\Inf$.
More precisely, there is a natural bijection $\varphi: V\Gamma \to \left\{  1, -1, \dots, k,-k \right\} \times [m]$ given by $\varphi\left( \varepsilon i, j
 \right) = ( \varepsilon j, i)$. 
 (with $\varepsilon \in \left\{ 1, -1 \right\}, i \in [m], j \in [k]$).  
 We thus need to find a suitable graph $\Sigma$ with vertex set $\left\{ 1, -1, \dots, k, -k \right\}$ such that the partition $\cP = \left\{ \left\{ j, -j \right\} : 1 \leq j \leq k \right\}$ is preserved by $\aut(\Sigma)$, and then show that $\varphi$ actually defines a graph isomorphism.

\begin{thm}\label{prop:graphs_inf}
 Let $\Gamma$ a vertex-transitive graph with vertex set $V\Gamma = \left\{ 1, -1, 2,-2, \dots,  m, -m \right\} \times [k]$ such that the group $G=\aut(\Gamma)$ satisfies \cref{eq:group_inf}. 
For $1 \leq j \leq k$, Denote by $D_{j}$ the set $\left\{ 1, -1, 2, -2, \dots, m,-m \right\} \times \left\{ j \right\}$ and let $\cD = \left\{ D_{j} : 1 \leq j \leq k \right\}$. 
  Then
  \begin{enumerate}
    \item\label{item:graphs_inf_block} If $1 \leq j_{1}, j_{2} \leq k$, then   $\Gamma(D_{j_{1}}) \cong \Gamma(D_{j_{2}})$ and $\Gamma(D_{j_{1}})$ is isomorphic to $mK_{2}$, $K_{m}\square K_{2}$, $\bar{mK_{2}}$ or $\bar{K_{m} \square K_{2}}$.
    \item \label{item:graphs_inf_graph} $\Gamma \cong \Inf[][][\Gamma/X^{k}][\cD/X^{k}]$ where
    \begin{equation}
      \label{eq:muKappa}
      (\lambda,\kappa)=
      \begin{cases}
        ( 1,0 )& \text{if } \Gamma(D)= mK_{2}\\
        ( 1,1 )& \text{if } \Gamma(D)=K_{m} \square K_{2} \\
        ( 0,1 )& \text{if } \Gamma(D)= \bar{mK_{2}}\\
        ( 0,0 )& \text{if } \Gamma(D)=\bar{K_{m} \square K_{2}}
      \end{cases}
    \end{equation}
    
    \item\label{item:graphs_inf_MinDeg4} The graph $\Gamma/X^{k}$ does not  admit a transposition preserving $\cD/X^{k}$.
  \end{enumerate}
\end{thm}
\begin{proof}
  Let us first proof \cref{item:graphs_inf_block}.
  It is clear that the two graphs $\Gamma(D_{j_{1}})$ and $\Gamma(D_{j_{2}})$ must be isomorphic, since $G$ acts transitively on $\cD$.
  Let $D \in \cD$. 
  Since $G_{D}^{D}=Y=\left\langle \sflip \right\rangle \times \sym(m) $, $\Gamma(D)$ must be isomorphic to one of the graphs listed in \cref{lem:graphs_easySF}.
  Assume that $\Gamma(D) \cong K_{m,m}$.
  Take two vertices $\alpha_{1}, \alpha_{2}$ in the same orbit of $X= G_{(V\Gamma \sm D)} = \left\{ 1 \right\} \times \sym(m)$. 
  Observe that for every $\beta \in D \sm \left\{ \alpha_{1}, \alpha_{2} \right\} $, $\alpha_{1} \sim \beta$ holds if and only if $\alpha_{2} \sim \beta$ holds.
  Moreover, since $G_{(\Omega \sm D)} \cong \left\{ 1 \right\} \times \sym(m)$, there exists an automorphism of $\Gamma$ that swaps $\alpha_{1}$ and $\alpha_{2}$ while fixing every vertex in $V\Gamma \sm D$.
  This implies that if $\gamma \in V\Gamma \sm D$, $\gamma \sim \alpha_{1}$ if and only if $\gamma \sim \alpha_{2}$.
  Therefore $\Gamma(\alpha_{1}) \sm \left\{ \alpha_{2} \right\} = \Gamma(\alpha_{2}) \sm \left\{ \alpha_{1} \right\}  $ and \cref{lem:graphsTrans} implies that $(\alpha_{1}\ \alpha_{2}) \in G$.
  This is obviously a contradiction to the fact that $G_{(\Omega \sm D)} =\left\{ 1 \right\} \times \sym(m)$.
  The same argument can be used to prove that $\Gamma(D)$ cannot be isomorphic to $K_{2m}$, $\bar{K_{m,m}}$ or $\bar{K_{2m}}$.

  Observe that the orbit of a vertex $(\varepsilon i, j) \in V\Gamma$ (with $\varepsilon \in \left\{ \pm 1 \right\}$, $1 \leq i  \leq m$ and $1 \leq j \leq k$) under the group $X^{k}$ is $\left\{ (\varepsilon \ell, j) : 1 \leq \ell \leq m \right\}$.
  This implies that each orbit of $\Gamma$ under $X^{k}$ (that is, a vertex of $\Gamma/X^{k}$) is determined by the pair $(\varepsilon, j)$ or equivalently, by the number $\varepsilon j$;
  let us denote such orbit by $A_{\varepsilon j}$.
  Observe that  $D_{j} = A_{j} \cup A_{-j}$ or equivalently, the partition $\cD/X^{k}$ of the vertices of $\Gamma/X^{k}$ is the family \[  \cD/X^{k} = \left\{ \left\{ A_{j}, A_{-j} \right\} : 1 \leq j \leq k \right\}  \]  
  Now consider the graph $\Upsilon = \Inf[][][\Gamma/X^{k}][\cD/X^{k}]$ where $(\lambda,\kappa)$ are as in \cref{eq:muKappa}.
  The vertices of $\Upsilon$ are the pairs $(A_{\varepsilon j }, i)$ with $i \in \left\{ 1, \dots, m \right\}$. 
  
  Now let us show that the mapping $\varphi: V\Gamma \to V\Upsilon$ given by $\varphi(\varepsilon i, j) = (A_{\varepsilon j }, i) $ defines an isomorphism.
  Assume that $v_{1}=(\varepsilon_{1} i_{1}, j_{1}) $ and $v_{2}=(\varepsilon_{2} i_{2}, j_{2})$ are such that $v_{1} \sim v_{2}$ in $\Gamma$. 
  If $j_{1} = j_{2}$, then the adjacency of $\varphi(v_{1})$ and $\varphi(v_{2})$ depends only on the pair $(\lambda,\kappa)$ but they were chosen precisely so that the graph $\Gamma[D_{j_{1}}]= \Gamma[D_{j_{2}}]$ is isomorphic to the graph $\Upsilon[\varphi(D_{j_{1}})] = \Upsilon\left[\left\{ A_{j_{1}}, A_{-j_{1}} \right\} \times [m]\right]$. 
  Assume thus that $j_{1} \neq j_{2}$ and observe that in this situation, $\left\{ A_{\varepsilon_{1} j_{1}}, A_{\varepsilon_{2} j_{2}}  \right\}\not\in \cD/X^{k}$ and $A_{\varepsilon_{1} j_{1}} \sim A_{\varepsilon_{2}j_{2}}$ in $\Gamma/X^{k}$, which implies that $\varphi(v_{1}) \sim \varphi(v_{2})$ in $\Upsilon$.

  The converse follows very similarly. 
  If $v_{1}=(\varepsilon_{1} i_{1}, j_{1}) $ and $v_{2}=(\varepsilon_{2} i_{2}, j_{2})$ are such that $\varphi(v_{1} ) \sim \varphi(v_{2})$ in $\Upsilon$ and $j_{1} \neq j_{2}$, then $A_{\varepsilon_{1} j_{1}} \sim A_{\varepsilon_{2} j _{2}}$ in $\Gamma/X^{k}$, which implies that there are $i'_{1}, i'_{2} \in \left\{1, \dots, m  \right\}$ such that $w_{1} = (\varepsilon_{1} i'_{1}, j_{1}) \sim (\varepsilon_{2} i'_{2}, j_{2}) = w_{2}$. 
  Finally notice that since $j_{1} \neq j_{2}$, there is an element of $X^{k}$ that maps $w_{1}$ to $v_{1}$ and $w_{2}$ to $v_{2}$, implying that $v_{1} \sim v_{2}$.
  If $j_{1} = j_{2}$, the fact that $v_{1} \sim v_{2}$ follows from the choice of $(\lambda, \kappa)$.

  Finally, \cref{item:graphs_inf_MinDeg4} follows from the fact that any transposition preserving $\cD/X^{k}$ must be a transposition swapping the two vertices of $\Gamma/X^{k}$ contained in a block of $\cD/X^{k}$. Those transpositions induce an automorphism of $\Gamma$ whose existence contradicts the fact that $G_{(V\Gamma\sm D)} = \left\{ 1 \right\} \times \sym(m)$.
\end{proof}

\cref{prop:graphs_inf} completes the classification of vertex transitive graphs whose automorphism group contains a group satisfying \cref{item:graphs4easy} of \cref{lem:graphs4} .
It just remains to classify the graphs whose automorphism group contains a group as in \cref{item:graphs4hard} of \cref{lem:graphs4}.
That is, we may assume that $\Gamma$ is a graph whose vertex-set $V\Gamma$ is $\left\{ 1, -1, 2, -2, \dots, m,-m \right\} \times [k]$ and $\aut(\Gamma)$ contains a group $G$ satisfying \[(\evenAb)^{k} \leq G \leq Y \wr \sym(k)\] for some $Y \leq \cubeWr$.
The following result classifies those graphs. 

\begin{prop} \label{prop:graphs_hard}
  Let $\Gamma$ be a vertex-transitive whose vertex set $V\Gamma$ is $\left\{ 1, -1, 2, -2, \dots, m,-m \right\} \times [k]$.
  Assume that $\aut(\Gamma)$ contains a group $G$ satisfying \[(\evenAb)^{k} \leq G \leq Y \wr \sym(k)\] for some $Y \leq \cubeWr$.
  If $m \geq 3$, then $\mu(\Gamma) = 2$.
\end{prop}
\begin{proof}
  Consider the vertices $\alpha_{1} = (1,1)$ and $\alpha_{2} = (-1,1)$.
  We shall prove that $\Gamma(\alpha_{1}) \sm \left\{ \alpha_{2} \right\} = \Gamma(\alpha_{2}) \sm \left\{ \alpha_{1} \right\}  $.
  Assume that $\omega \in \Gamma(\alpha_{1}) \sm \left\{ \alpha_{2} \right\} $.
Since $m \geq 3$, there must exist $i \in \left\{ 2, \dots, m \right\}$
  such that  $\omega \not\in \left\{ (i,1), (-i,1) \right\}$.
  Let  $\gamma_{1}= (i,1)$ and $\gamma_{2} = (-i,1)$.
  By our assumptions on  $G$,  the permutation $(\alpha_{1} \ \alpha_{2})(\gamma_{1}\ \gamma_{2})$ is an element in $G$.
  This permutation fixes $\omega$ while swapping $\alpha_{1}$ and $\alpha_{2}$ proving that $\omega \in \Gamma(\alpha_{2}) \sm \left\{ \alpha_{1} \right\} $.
  The other inclusion is analogous.
\end{proof}

Finally, observe that if $\Gamma$ satisfies the hypothesis of \cref{prop:graphs_hard} with $m=2$ then, up to a relabelling of the vertices, $\Gamma$ satisfies the assumptions of \cref{prop:graphs_inf}. 
It is enough to observe that groups $\evenAb[2]$ and $\left\{  1\right\} \times \sym(2)$ acting on $\left\{ 1,-1 \right\} \times \left\{ 1,2 \right\}$ are permutation isomorphic.
This concludes the proof of \cref{thm:graphs}.

\printbibliography

\end{document}